\documentclass{article}
\usepackage{plain}

\begin{document}

\newcommand{\n}{\noindent}
\newcommand{\cl}{\centerline}
\newcommand{\vp}{\varepsilon}
\newcommand{\eightrm}{}
\newcommand{\ninerm}{}

\begin{plain}
%%  \input{plain-doc}
%\end{plain}

%mart-final-backup-Mar25.tex

\parskip=6pt
\hbadness=10000
\vbadness=10000
\baselineskip = 18pt
\def\n{\noindent}
\def\vp{\varepsilon}
\overfullrule = 0pt

%%%%\runningauthor={\eightrm 

\font\ninerm=cmr10 at 9 pt

\font\eightrm=cmr10 at 8 pt

%%%%{\cmcsc10 Ro}

\outer\def\proclaim  #1 #2 #3. #4\par{\medbreak
 %%%\noindent
 #1{\eightrm #2} {\rm #3.}\enspace{\sl#4}\par
 \ifdim\lastskip<\medskipamount \removelastskip\penalty55\medskip\fi}

\def\\X {\cal \X }
\def\cl{\centerline}

\def\skp{\bigskip}
\def\sqr{\quad  \vrule height 6pt depth 0pt width 6pt}

\def\E{{\rm {\bf E}}}
\def\Z{{\rm {\bf Z}}}
\def\q{\  }
\def\l{\leqalignno}
\def\({\left(}  \def\){\right)}   \def\e{{\rm e}}
\def\Var{{\rm Var}}  \def\P{{\rm {\bf Pr\,}}} \def\R{{\rm {\bf R}}}
\def\D{{\rm {\bf D}}}
\def\Re{\R}
\def\X {{\cal X }}
\def\\X {{\cal X }}
\def\t{\tilde}

\def\i#1 {\item{\hbox to 1cm {#1\hfil}   }  }

\def\a{\alpha}
\def\<{\preceq}  \def\>{\succeq}   
\def\e{{\rm e}}
\def\len{\left\|}
\def\rn{\right\|}

\def\sinf{ \vphantom{ \sup     } \inf    }
  \def\Y#1  {\vphantom { {\scriptstyle Y_2  }           }   #1 }
  \def\nsup{\sup\nolimits}      \def\ninf{\inf\nolimits}
\def\sinf{ \vphantom{ \sup     } \inf    }
  \def\Y#1  {\vphantom { {\scriptstyle Y_2  }           }   #1 }
  \def\nsup{\sup\nolimits}      \def\ninf{\inf\nolimits}
\def\vpi{\varphi}
\def\d{{\rm d}}

\cl {\bf OPTIMUM BOUNDS FOR THE DISTRIBUTIONS OF MARTINGALES}
\cl {\bf  IN BANACH SPACES$^1$}

\vskip .5cm

\cl{B{\eightrm Y}  I{\eightrm OSIF} P{\eightrm INELIS}} \footnote
{}{\parskip 1pt {\ninerm ${}^1$Research partially supported by NSF grant  DMS-8921369.


moments, martingales in  Banach spaces, 2-smooth  Banach spaces, sums of independent random variables.}}%END FOOTNOTE

\cl{{\it Michigan Technological University}}

\font\ninerm=cmr10 at 9 pt

%%%%%%%%%%%\font\twicemathfont=mathfont at 15pt

{\ninerm

 A general device is proposed, which provides for extension of exponential inequalities for sums of independent real-valued random variables to those for martingales in the 2-smooth Banach spaces. This is used to obtain optimum bounds of the Rosenthal-Burkholder and Chung types on moments of the martingales in the 2-smooth Banach spaces. In turn, it leads to best-order bounds on moments of the sums of independent random vectors in any separable Banach spaces. Although the emphasis is put on the infinite-dimensional martingales, most of the results seem to be new even for the one-dimensional ones. Moreover, the bounds on the Rosenthal-Burkholder type of moments seem to be to certain extent new even for the sums of independent real-valued random variables. Analogous inequalities for (one-dimensional) supermartingales are given.} 

\vskip .7cm

\n {\bf 1. Introduction.}
\ \  For a separable Banach space $(\X ,\ \|\cdot\|)$, let ${\cal S}(\X )$ denote
the
class of all 
%%%%``stopped'' 
 sequences $f=(f_j) = (f_0,f_1,\ldots )$
 of Bochner-integrable random vectors in $\X $, with
$f_0\equiv 0$, defined on a probability space $(\Omega, F, P)$ and adapted to a
non-decreasing sequence $(F_j) = (F_0,F_1,\ldots)$  of
sub-$\sigma$-fields
of $F$; here, $n$, $(\Omega, F,P)$ and $(F_j)$ are considered
attributes of
$f$ and may be different for different $f\in {\cal S}(\X )$. 
%%%We also
%%%put
%%%$f_\infty = f_n$ for $f\in {\cal S}(\X )$.

%%We restrict the forthcoming consideration to that of the stopped
%%%sequences
%%with the only purpose of simplifying the wording. Also, the results
%%can be
%%easily extended to any measurable seminormed spaces $\X $.

For $f\in {\cal S}(\X )$, put $f^* = \sup\{\|f_j\|\colon \ j=0,1,\ldots\}$,
$d_0=
d_0(f) \equiv 0$, $d_j = d_j(f) = %\break 
f_j-f_{j-1}$,
$j=1,2,\ldots$,
$S_p =
S_p (f) = \left(\sum ^\infty_{j=1} \|d_j\|^p\right)^{1/p}$,
$p>0$,
$s_2=s_2(f)
\hfill\break = \left(\sum ^\infty_{j=1} \E _{j-1}\|d_j\|^2\right)^{1/2}$,
where
$\E _{j-1}$ stands for the conditional expectation given $F_{j-1}$.

Let ${\cal {\cal M}  }(\X )$ denote the class of all sequences $(f_j) \in
{\cal
S}(\X )$ that are martingales  and ${\cal {\cal M}  }_{ind}(\X )$,  the class
of all
sequences $(f_j) \in {\cal S}(\X )$ having independent increments
$d_j$'s.

%\headline={\ifnum\pageno>1\ifodd\pageno
%\rightheadline \else\leftheadline\fi\else{}\fi}
\def\rightheadline
{\ninerm\hfil INEQUALITIES FOR VECTOR MARTINGALES
\hfil}
\def\leftheadline{\ninerm\hfil I. PINELIS \hfil}
\voffset=-.5\baselineskip

For any two nonnegative expressions ${\cal E}_1 $ and ${\cal E}_2$,
let us
write ${\cal E}_1\< {\cal E}_2$ (or, equivalently, ${\cal E}_2\> {\cal
E}_1$)
if ${\cal E}_1\le A{\cal E}_2$, and ${\cal E}_1\asymp {\cal E}_2$ if
${\cal
E}_1 \< {\cal E}_2 \< {\cal E}_1$. Here, $A$
denotes
a positive absolute constant.

%%%In this paper, when we speak about optimum bounds, we mean the
%%%optimality
%%%%%to the equivalence $\asymp$.

We assume that $\inf \emptyset =\infty$, $\sup\emptyset=0$,
$\sum_{j\in\emptyset} u_j = 0$, $\prod_{j\in
\emptyset}
u_j=1$.

Symbol \sqr \ \ indicates the end of a proof or a remark.

In section 2, some preliminary results on  the 2-smooth Banach spaces and on the martingales in such spaces are given.

In section 3,  a device is suggested, which provides for the extension of
exponential inequalities for the sums of independent real-valued
random
variables to those for the martingales in the 2-smooth Banach spaces. In
particular, by that means an exponential inequality for the
martingales in
the 2-smooth spaces optimal in terms of $\|d^*\|  _\infty$,
$\|  s_2\|  _\infty$ is obtained, which is a generalization of an
inequality
of Bennett (1962) and Hoeffding (1963).

In section 4, using methods of Burkholder (1973) and results of
section~3,
 we obtain optimal (to the above defined relation $\asymp$) upper
bounds
of the Rosenthal (1970)-Burkholder (1973) type on moments of the
martingales in the 2-smooth Banach spaces, i.e., optimal in terms of
$\|  d^*\|  _p$, $\|  s_2\|  _p$, for $p\ge 2$.

In section 5, via a modification of the martingale decomposition
method by
Yurinskii (1974), we apply the inequalities of section~4 to obtain
bounds
of the Rosenthal-Burkholder type on $\|  \ \|f_n\|-\E\|f_n\|\ \|  _p$,
$p\ge
2$, for an arbitrary separable Banach space $\X$, but only for $f\in
{\cal M}  _{ind}(\X )$.

In section 6, we show that the inequalities derived in sections~4, 5
are
optimal in the terms used. The key roles in the proof of the
optimality are
played by the results and methods of Pinelis and Utev (1984).
%% the sums of independent real-valued random variables with a common symmetric
%%distribution concentrated at three points turn out to represent the case
%%when the moments    are in certain sense largest.

In section 7, we obtain bounds on $\| \, \|  f_n\| \, \| _p$, $p\ge 2$, $f\in
{\cal M}  (\X )$,
which are optimal in terms of $n$ and $\|S_p\|_p$.  We refer to them as to
bounds
of  Chung type.

In section 8, inequalities for super-martingales (of course, in
$\X =\Re
%{\Bbb R}
$) similar to those in sections 3, 4 and  certain refinements for real-valued martingales are presented.

In section 9, some details are given.	For instance, it is shown
that, at
least for the sums of independent random vectors,
the
bounds in Pinelis (1980) are equivalent to those in section~4. 

Substantial part of the results was announced in Pinelis (1992).

I am grateful to Donald Burkholder and Walter Philipp for their
support and
stimulating conversations and to Pawel Hitczenko, who has acquainted
me
with his latest results here referred to. I would like to thank the organizers of the work-shop "Linear Analysis and Probability" in the Department of Mathematics, Texas A\&M University for the great opportunity to discuss the results of this paper.

%%%%%%%%%%%%%%%%%%%%%%%%%%%%%%%%%%%%%%%%%%%%%%\bye
\vskip .8cm

\n {\bf 2.  Preliminaries: 2-smooth Banach spaces and a reduction of martingales.}
\ \ Let us call a Banach space $(\X , \, \|\cdot\|)$ 2-smooth [or, more exactly, $(2,D)$-smooth, for some $D = D(\X ) >0$] if $\ \forall x\in \X \ \forall y\in \X $\ \  
$$\l{&\| x+y\| ^2 + \| x-y\| ^2 \le 2\| x\| ^2 +2D^2\| y\| ^2. &(2.1)}$$

The importance  of the 2-smooth spaces was elucidated in the paper of Pisier (1975): they play the same role with respect to the vector martingales  as the spaces of type 2 do with respect to the sums of independent random vectors.

The definition assumed in this paper is slightly different from that given by Pisier [which required only that (2.1) hold for an equivalent norm], because, in the subsequent account, we would like to follow the dependence of certain constants  on $D$, the constant of the 2-smoothness.

Substituting $\lambda x$ for $y$ in (2.1), where $\lambda \in \Re$, one observes  that except for $\X=\{0\}$, $D\ge 1$.

It is easily seen that the condition 
$$ \l{&(\| x\| ^2)''(v,v) \le 2D^2 \| v \| ^2 \quad \forall x\in \X \ \forall v\in \X &(2.2)}$$
is sufficient for the $(2,D)$-smoothness, where $(\| x\| ^2)''(v,v) $ stands for the second directional derivative of the function $x \mapsto \| x\| ^2$ in the direction $v$. 

By way of illustration, we give

\proclaim P ROPOSITION 2.1.  For any $p\ge 2$ and any measure space
$(T,{\cal A}, \nu)$, $L^p := L^p (T, {\cal A}, \nu)$ is $(2,\sqrt {p-1} \, )$-smooth.

P{\eightrm ROOF.} For  $\|x\| := (\int_T |x|^p d\nu)^{1/p}$,   one has
$$\leqalignno{
{1\over 2} (\| x\| ^2)''(v,v)&=
(p-1)\|x\|^{2-p} \int_T |x|^{p-2} v^2\d\nu\cr
&\quad - (p-2) \|x\|^{2-2p}
\left(\int_T|x|^{p-2}xv\,\d\nu\right)^2\cr &\le (p-1)\|v\|^2-0&(2.3)}$$
if $x\in \X \backslash \{0\}$,  $v\in \X $ (
in view of
H\"older's inequality), and ${1\over 2} (\| x\| ^2)''(v,v) =\|v\|^2$ if $x=0$.  Thus, (2.2) is checked. \sqr

In particular, what is obvious and well-known, if $\X $ is a Hilbert space, then it is (2,1)-smooth.

\n R{\eightrm EMARK. } If $L^p = L^p (T, {\cal A}, \nu)$ is at least
two-dimensional, i.e., if there exist $T_1\in{\cal A}$, $T_2\in{\cal A}$ such that $T_1 \cap T_2 = \emptyset$, $0<\nu (T_1) <\infty $, $0<\nu (T_2) <\infty $, 
 then $L^p$ is {\it not \/} $(2,D)$-smooth if $0< D< \sqrt{p-1}$,  so that Proposition~2.1 gives the best bound.
Indeed,
put $x\equiv \nu (T_2)^{1/p}$ on $T_1$, $x\equiv
\nu(T_1)^{1/p}$ on $T_2$, $v\equiv x$ on $T_1$, $v\equiv
-x$ on
$T_2$, $x\equiv v \equiv 0$ on $T\backslash T_1 \backslash T_2$. Then $x\ne 0$, and (2.3) turns into an equality. \sqr

Condition (2.2) is not only sufficient but also  necessary for the $(2,D)$-smoothness if the derivatives are understood in a generalized sense. To state this remark rigorously, let us give a definition somewhat extending the notion of the 2-smoothness.

For  any  Banach space $(\X , \|\cdot\|)$, we call a function $\Psi \colon \X    \to [0,\infty)$ $(2, D)$-smooth, $D>0$, if it satisfies the conditions: $\Psi (0)=0$,
$$ |\Psi (x +v) -\Psi (x )| \le \| v \| ,$$
$$ \Psi  ^2 (x +v) - 2 \Psi ^2 (x ) + \Psi ^2 (x -v)    \le 2 D ^2 \| v \| ^2$$
for all $x \in \\X   $, $v\in \\X   $. 

Evidently, a Banach space $(\X , \|\cdot\|)$ is $(2, D)$-smooth if and only if its norm function  is $(2, D)$-smooth.

The results stated in the subsequent sections for the norms of the martingales in the 2-smooth spaces can be  extended to those for the $(2, D)$-smooth functions of martingales in any Banach spaces.

For any $(2, D)$-smooth function $\Psi$ on a {\it finite-dimensional\/} Banach space $(\X , \|\cdot\|)$ and $\vp >0$, define 
$$\l{&\Psi _\vp (x) = \sqrt{ \int _{\X}  \Psi ^2 (x -\vp y) \gamma (\d y)  },   &}$$
where  $\gamma $ is, say, a zero-mean Gaussian measure on $\X$ with ${\rm support} (\gamma) = \X$.

\proclaim L EMMA 2.2. If  $ \Psi $ is a $(2, D)$-smooth function on a finite-dimensional Banach space \break $(\X , \|\cdot\|)$, then for all $\vp >0$, $ \Psi _\vp$ has the Fr\'echet 
derivatives $\Psi _\vp '(x)$, $\Psi _\vp '' (x), \ldots$ of any order, and 
the directional derivatives in any direction $v \in \X $ satisfy the inequalities 
$$\l{&|\Psi _\vp '(x)(v)| \le \| v \|, \qquad (\Psi _\vp ^2 )''(x)(v,v) \le 2 D ^2 \| v \| ^2  &(2.4)}$$
 for all $x\in \X$. Besides, for each $x\in \X$, $\Psi _\vp (x) \to \Psi  (x) $ as $\vp \downarrow 0$. [In this generalized sense,  sufficient condition (2.2) is also necessary for a Banach space $(\X , \|\cdot\|)$ to be $(2, D)$-smooth; note that (2.1) may be considered locally -- for any two-dimensional subspace containing $x$ and $y$.]

P{\eightrm ROOF. } Among the statements of the lemma, only the first of the inequalities (2.4) is comparatively non-trivial. Observe that 
$$\l{  |(\Psi _\vp ^2 )'(x)(v)| &\le \mathop { {\rm lim }\, {\rm sup } }\limits _{t\to 0} \; {1\over t}  \int _\X |\Psi   (x+tv-\vp y)-\Psi   (x-\vp y)| \(\Psi   (x+tv-\vp y)+\Psi   (x-\vp y)\) \gamma (\d y)
 \cr
& \le  \| v \| \int _\X 2 \Psi   (x-\vp y)  \gamma (\d y). \cr }$$
One can assume that $\Psi _\vp  (x) \ne 0$ for all $x\in \X$. Now, 
$$\l{  |(\Psi _\vp  )'(x)(v)| &= \left| {   {(\Psi _\vp ^2 )'(x)(v)} \over 
{2\Psi _\vp  (x)}  } \right| \cr
&\le {  { \| v \| \int _\X 2 \Psi   (x-\vp y)  \gamma (\d y)    }\over { 2 \sqrt { \int _{\X}  \Psi ^2 (x -\vp y) \gamma (\d y)}   }   }\cr
&\le \| v \| \cr}$$
by Schwartz inequality. \sqr

We also need the following folklorish lemma.

\proclaim L EMMA 2.3. Let $(f_j)_{j=0}^\infty \in {\cal M}  (\X )$ be a martingale in a separable Banach space  $(\X , \|\cdot\|)$ relative to a filtration $(F_j)_{j=0}^\infty$. Then for any $\vp >0$, there exists a martingale $(f_{j,\vp})_{j=0}^\infty \in {\cal M}  (\X )$ relative to a filtration $(F_{j,\vp})_{j=0}^\infty$ such that $\forall j=1,2, \ldots , $ 
$$\l{& f_{j,\vp}\ {\rm is\ a\  random\ variable\ having\ only\ a\ finite\ number\ of\ values,}&(2.5)\cr
 &f_{j,\vp} \to f_j\ {\rm in\ probability\ as\ } \vp \downarrow 0,&(2.6)\cr
 &\E g(f_{j,\vp}) \le \E g(f_{j}), &(2.7)\cr
 & \|\E (g(f_{j,\vp})|F_{j-1,\vp}) \| _\infty \le \E (g(f_{j})|F_{j-1}) \| _\infty , &(2.8)\cr }$$
 where $g$ is any nonnegative convex real function on $\X$.

P{\eightrm ROOF. } Consider the approximation $f_{j,\vp} :=  \E (f_j | F_{j,\vp})$, where $F_{j,\vp}$ is the $\sigma$-field generated by all the events of the form $\{  f_i \in B_{k,\vp} \} $, $i=0,1, \ldots , j$, $k=1, \ldots , k(j,\vp)$,  where  $(\{ B_{k,\vp} \colon k=1, \ldots , k(j,\vp )\})_{j=0}^\infty$ is an increasing sequence of sets of balls in $\\X   $ of the radius $\vp$ such that $k(0,\vp) < k(1,\vp) <  k(2,\vp) < \ldots $ and
%%\hfill 
 $\P ( \bigcap _{i=1}^j \bigcup _{k=1}^{k(j,\vp)} \{ f_i \in B_{k,\vp} \} ) \ge 1-\vp$, $j=1,2, \ldots $ [the existence of such a sequence of sets is guaranteed by the tightness of any probability measure on a separable Banach space].
 
 Then, (2.5) and (2.6) are satisfied. The Jensen's type inequality
 $g(f_{j,\vp}) \le \E (g(f_{j})|F_{j,\vp})$
 implies (2.7) and 
 $\| \E (g(f_{j,\vp})|F_{j-1,\vp}) \| _\infty \le \| \E (g(f_{j})|F_{j-1, \vp}) \| _\infty \le \| \E (g(f_{j})|F_{j-1}) \| _\infty ,$ 
 so that (2.8) is also true. \sqr 

\n R{\eightrm EMARK} 2.4.  Using the standard construction
$(f_{j\land n \land \tau })_{j=0}^\infty$, where $n=1, 2, \ldots  $, \break $\tau = \inf \{ j\colon \|f_j\| \ge M \}$ with large $n$ and $M$, we will need only to consider the bounded stopped martingales. By Lemma 2.3 and Lemma 2.2, when proving   subsequent results for the martingales in the 2-smooth spaces, we may and shall, without loss of generality, restrict the consideration only to the bounded stopped martingales $(0,f_1, \ldots , f_n, f_n, f_n, \ldots )$ with each of $f_j$'s having only a finite number of values in a finite-dimensional Banach space $(\X , \, \|\cdot\|)$ satisfying condition (2.2); for such martingales, we put $f_\infty =f_n$. \sqr

%%%%%& F_{j,\vp} {\rm is a sub-$\sigma$-algebra of F_j},&(2.6)\cr

\proclaim P ROPOSITION 2.5. If $f\in {\cal M}  (\X )$ and $\X $ is a $(2,D)$-smooth
separable
Banach space, then $$\|  f_n\|  _2 \le D\|  s_2\|  _2 = D\|  S_2\|  _2.$$

P{\eightrm ROOF.}  For $j=1,2,\ldots$, put $$g(t) = \E \|f_{j-1} +
td_j\|^2.$$
In view of Remark 2.4, we may assume that $g'(0) = 0, g''(t) \le 2D^2\|  d_j\|  ^2_2$, and so,
$$\eqalign{\E \|f_j\|^2 - \E \|f_{j-1}\|^2 &= g(1) - g(0) = \int^1_0
g''(t)
(1-t)dt\cr &\le D^2\|  d_j\| _2  ^2, \quad j=1,2,\ldots\ .}$$ It remains
to sum
up these inequalities. \sqr

\medskip

The upper bounds provided by Theorems~15.1 and 21.1 of Burkholder
(1973)
can be immediately extended to the martingales in the 2-smooth
separable
Banach spaces.	Let us state this as

\proclaim T HEOREM 2.6. 	If $f\in {\cal M}  (\X )$, $\X $ is a 2-smooth separable
Banach
space, a function $g \colon \ [0,\infty) \to [0,\infty)$ is
non-decreasing, $g (0) = 0$, and $g (2u) \le
c_g g (u)$,
$u\ge 0$, for some $c_g $, then $$\E g (f^*) \le
c_{2.9}[\E g (D\,s_2) + \E g (d^*)];\leqno (2.9)$$ if, moreover,
$g $ is convex, then
$$\E g (f^*) \le c_{2.10} \E g (D\,S_2);\leqno
(2.10)$$ here, $c_{2.9}$ and $c_{2.10}$ depend only on $c_g $.

P{\eightrm ROOF}  repeats that in Burkholder (1973) with the following
exceptions: 1) use
$\|\cdot\|$ instead of $|\cdot|$, and 2) use Proposition~2.5 instead of the identities (in the notation
therein):\ $\|h\|_2 = \|S(h)\|_2$ and $\|h\|_2 = \|s(h)\|_2$. \sqr

It is well-known [see (11.1) in Burkholder (1991)] that
$(Ap)^p$ is
optimum for $c_{2.10}$ at least if $\X $ is a Hilbert space.  As to
optimum
bounds like (2.10) for conditionally symmetric martingales, see
section~4 below. 
 
\vskip .8cm

\n {\bf 3.  Exponential bounds on tail probabilities for the martingales
in the
2-smooth spaces.}

\proclaim T HEOREM 3.1.	Suppose that $f\in {\cal M}  (\\X  )$, $\\X  $ is a $(2,D)$-smooth separable
Banach space, and $\lambda>0$ is such that
$\E\e ^{\lambda\|d_j\|} <
\infty$ for $j=1,2,\ldots$\ .  Then for all $r\ge 0$,
$$\leqalignno{\P (f^*\ge r) &\le 2\e ^{-\lambda r} \left|\left|
\prod
^\infty_{j=1} [1+D^2\E _{j-1} (\e ^{\lambda\|d_j\|} - 1 -
\lambda\|d_j\|)]\right|\right|_\infty& \cr &\le 2
\exp\left\{-\lambda r
+ D^2 \left|\left| \sum^\infty_{j=1} \E _{j-1} (\e ^{\lambda
\|d_j\|} - 1
- \lambda\|d_j\|)\right|\right|_\infty\right\}.& }$$

We shall obtain this theorem as a particular case of the following result for the \break $(2,D)$-smooth functions defined in section 2.

\proclaim T HEOREM 3.2 .	Suppose that $f\in {\cal M}  (\X   )$, $\\X   $ is any separable
Banach space, $\lambda>0$ is such that
$\E\e ^{\lambda\|d_j\|} <
\infty$ for $j=1,2,\ldots$\ , and a function $\Psi $ is $(2,D)$-smooth. 
Then for all $r\ge 0$,
$$\leqalignno{\P \(\sup_j \Psi (f_j) \ge r\) &\le 2\e ^{-\lambda r} \left|\left|
\prod
^\infty_{j=1} (1+e _j)\right|\right|_\infty, &(3.1) \cr
 &\le 2
\exp\left\{-\lambda r
+ \left|\left| \sum^\infty_{j=1} e _j \right|\right|_\infty\right\},&(3.2) }$$
where 
$$ e _j =  {D_*^2 \,  } \E _{j-1} (\e ^{\lambda \|d_j\|} - 1 -
\lambda  \|d_j\|),\quad D_* = 1 \vee D.$$

   P{\eightrm ROOF. } Put $u(t)=u_{x ,v} (t)=\Psi (x +tv)$ for any $x ,v$ in $\\X   $. By Lemma 2.2 and Remark~2.4, we may assume that $u$ is Fr\'echet-differentiable, 
    $|u'(t)|\le \| v\|, \ \; \break (u^2)''(t)\le 2 D^2 \| v\| ^2.$
Hence, when $u''u>0$, one has  
$(\cosh u)'' =u'{}^2
\cosh u + u'' \sinh u \break \le (u'{}^2+ u''u)\cosh u 
={1 \over 2} (u^2)'' \cosh u$
%%\break 
$\le  D^2 \| v\| ^2 \cosh u , 
$  
and, otherwise,
$(\cosh u)'' \le\break u'{}^2
\cosh u \le \| v\| ^2 \cosh u ; $
so, in any case,
$$\l{&(\cosh u)''\le D_*^2 \,  \| v\| ^2 \cosh u. &(3.3)\cr }$$ 
  Consider now
$\varphi(t) := \E _{j-1} \cosh(\lambda \Psi (f_{j-1} + td_j )), \quad
|t|\le 1.$ 
 In view of (3.3) and Remark 2.4,
$$\eqalign{\varphi''(t) &\le D_*^2 \, \lambda^2 \E _{j-1} \|d_j\|^2
\cosh (\lambda\Psi (f_{j-1} + td_j))\cr &\le D_*^2 \,  \lambda^2 \E _{j-1}
\|d_j\|^2
\e ^{\lambda  t\|d_j\|} \cosh (\lambda \Psi (f_{j-1})),}$$
$|t|<1$.  But $\varphi'(0)
= 0$ since $(f_j)$ is a martingale, and therefore,
$$\E _{j-1} \cosh (\lambda\Psi (f_j)) = \varphi(1)
=\varphi(0) + \int^1_0 (1-t) \varphi''(t)\d t \le (1+e _j) \cosh (\lambda
\Psi (f_{j-1})).$$

Thus, putting $G_0=1$,
$G_j =
\cosh (\lambda\Psi (f_j))/ \prod^j_{i=1} (1+e _i),$ 
$j=1,2,\ldots,$ one has a positive supermartingale.  Hence, if
$\tau := \inf\{j\colon \
\Psi (f_j)\ge r\},$ 
then $\E G_\tau \le \E G_0 = 1$, and so,
$$ \P \(\sup_j \Psi (f_j) \ge r\) \le 
  \P\left(G_\tau \ge \cosh (\lambda r)\left/ \left|\left|\prod
^\infty_{j=1} (1+e _j)\right|\right|_\infty\right)\right.  . 
$$
 Now, (3.1) follows from Chebyshev's inequality and $\cosh u >
\e ^u/2$; (3.2)  is elementary. \sqr

%%%\bye%%%%%%

\n R{\eightrm EMARK.} For $f\in {\cal M}  _{ind}(\Re
%{\Bbb R}
)$, i.e., for the sums of
independent zero-mean real-valued random variables, the following is
used as
a starting point when proving exponential inequalities:
$$\eqalign{\P (f_n\ge r) &\le \e ^{-\lambda r} \prod^n_{j=1}
[1+\E (\e ^{\lambda
d_j} - 1 - \lambda d_j)]\cr &\le \exp \left\{-\lambda r +
\sum^n_{j=1}
\E (\e ^{\lambda d_j} - 1 - \lambda d_j)\right\}.}$$

Thus, Theorem 3.1 provides a similar starting point for $f\in {\cal M}  (\X )$,
$\X $
being \break 2-smooth.  [In this sense, it is analogous to the results of
Pinelis
and Sakhanenko (1985) for the sums of independent random vectors.] A
general
method of obtaining exact exponential inequalities for the sums of
independent
real-valued random variables is proposed in Pinelis and Utev (1989).
So,
for the martingales in the 2-smooth spaces, these two devices taken
together produce analogues of exact ``independent real-valued''
exponential
bounds. For instance, this remark easily leads to the following
analogues
of classical results of Bernstein 
%%%(1924)
 and Bennett (1962)-Hoeffding (1963) [cf. Theorem~9 and 3, resp., in Pinelis and Utev (1989)]. \sqr

\proclaim T HEOREM 3.3.	Suppose that $f\in {\cal M}  (\X )$, $\X $ is a $(2,D)$-smooth
separable Banach space, and $$\left\| \sum^\infty_{j=1} \E _{j-1}
\|d_j\|^m\right\|_\infty \le m!	\Gamma^{m-2}B^2/(2D^2)$$ for
some
$\Gamma>0$, $B>0$ and $m=2,3,\ldots\ $.  Then for all $r\ge 0$,
$$\P (f^*\ge
r) \le 2\exp \left( - {r^2 \over B^2+B\sqrt{B^2+2\Gamma r}}\right).$$

P{\eightrm ROOF.} Under the conditions given, 
$$D^2\left\|
\sum^\infty_{j=1} \E _{j-1} (\e ^{\lambda\|d_j\|} - 1- \lambda
\|d_j\|)\right\|_\infty \le {1\over 2} \sum^\infty_{m=2}
\lambda^m
\Gamma^{m-2}B^2 ={B^2\lambda^2\over 2(1-\lambda\Gamma)}, 0 \le \lambda
<{1\over \Gamma}.$$ 
Now Theorem 3.1 yields $$\P (f^*\ge r) \le 2 \exp\left\{
-\lambda
r + {B^2\lambda^2 \over 2(1-\lambda\Gamma)}\right\}.$$ It remains to
minimize the R.H.S.  in $\lambda$. \sqr

\proclaim T HEOREM 3.4.	Suppose that $f\in {\cal M}  (\X )$, $\X $ is a $(2,D)$-smooth
separable Banach space, and $\|  d^*\|  _\infty \le a$,
$\|  s_2\|  _\infty \le
b/D$ for some $a>0$, $b>0$.  Then for all $r\ge 0$,
$$\leqalignno{\P (f^*\ge
r) &\le 2 \exp\left[ {r\over a} - \left({r\over a} + {b^2\over
a^2}\right)
\ln \left(1+{ra\over b^2}\right)\right]&(3.4)\cr &\le 2
\left({\e b^2\over
ra}\right)^{r/a}.&(3.5)}$$

P{\eightrm ROOF.} Since the function $g(u) := u^{-2}(\e ^u-1-u)$ for
$u\ne 0$,
$g(0) := 1/2$ is increasing in $u \in \Re
%{\Bbb R}
$, $$\E _{j-1}
(\e ^{\lambda\|d_j\|} - 1 - \lambda\|d_j\|) \le {\e ^{\lambda a} - 1 -
\lambda
a\over a^2} \E _{j-1} \|d_j\|^2.$$ Now Theorem 3.1  yields $$\P (f^*\ge r) \le 2
\exp
\left\{ - \lambda r + {\e ^{\lambda a}-1-\lambda a\over a^2}
b^2\right\},$$
and the minimization in $\lambda$ gives (3.4).  Inequality (3.5) is
trivial.
\sqr

In the special case $\X =L^2$, a bound similar to (3.4), but somewhat weaker, was proved by
Kallenberg
and Sztencel (1991);  their method seems
to be
confined only to Hilbert spaces.

Theorem 3.4 was proved in Pinelis (1992) for $\X =L^p$, $p\ge 2$; a
version
for the general 2-smooth spaces was given therein too, but with
another,
greater constant in place of $D$.

\proclaim T HEOREM 3.5. 	Suppose that $f\in {\cal M}  (\X )$, $\X $ is a $(2,D)$-smooth
separable Banach space, and $\sum\limits^\infty_{j=1}
\|  d_j\|  ^2_\infty
\le b^2_*$ for some $b_* >0$.  Then for all $r\ge 0$, $$\P (f^*\ge r) \le 2 \exp\left\{ -
{r^2
\over 2D^2b^2_*}\right\}.$$

P{\eightrm ROOF}  is the same as that of Theorem 3 in Pinelis (1992) except that,
in
view of (3.3), one can use $D^2$ instead of $B$ therein. \sqr

Theorem 3.5 can be improved in the special case of conditionally
symmetric
martingales.

\proclaim T HEOREM 3.6.	Suppose that $\X $ is a $(2,D)$-smooth separable
Banach
space, $f\in {\cal M}  (\X )$, $\|  S_2(f)\|  _\infty \le b$ for some $b >0$, and the increments
$d_j$'s
are $F_{j-1}$-conditionally symmetrically distributed,
$j=1,2,\ldots$\ .
Then for all $r\ge 0$, $$\P (f^*\ge r) \le 2 \exp\left\{ - {r^2\over
2D^2b^2}\right\}.$$

P{\eightrm ROOF.} Being conditionally symmetric, $(f_j)$ is a
martingale also
relative to the sequence $(G_j)$, where $G_j$ is the $\sigma$-field
generated by $F_j$ and $\|d_{j+1}\|$; see, e.g.
Lemma~10.2
in Burkholder (1991).  Now, the proof can be concluded as that of
Theorem~3 in Pinelis (1992):\ only the conditional expectations
given
$G_j$'s are taken instead of those given $F_j$'s, and $D^2$,
$\|d_j\|^2$
are used in place of $B$, $b^2_n$ therein, resp. \sqr

In the case when $\X =\Re
%{\Bbb R}
$ and $d_j$'s are simple functions,
Theorem~3.6
was given in Hitczenko (1990a). 
%%%%Upper bounds for the $L_p$-norms of
%%%martingales, Probab.  Th.  Rel.  Fields  {\bf 86}, 225-237.

An analogous result for sums of independent random vectors in
arbitrary
separable Banach spaces is given in Pinelis (1990). 

\vskip .8cm

\n {\bf 4. A spectrum   of Rosenthal-Burkholder type bounds on
moments of the 
martingales in the 2-smooth spaces.}

\proclaim T HEOREM 4.1.	If $f\in {\cal M}  (\X )$, $\X $ is a $(2,D)$-smooth separable
Banach
space, $p\ge 2$, $1\le c\le p$, then
$$\|  f^*\|  _p \< c\|  d^*\|  _p + \sqrt
c\, \e ^{p/c} D\|  s_2\|  _p;\leqno (4.1)$$
 in particular,
$$\leqalignno{\|  f^*\|  _p &\< p\|  d^*\|  _p + \sqrt p\,
D\|  s_2\|  _p,&(4.2)\cr \|  f^*\|  _p &\< {p \over \ln p} (\|  d^*\|  _p +
D\|  s_2\|  _p),&(4.3)\cr \|  f^*\|  _p &\< \alpha\|  d^*\|  _p +
\e ^{p/\alpha}D\|  s_2\|  _p, \quad 1\le \alpha \le {p\over \ln (\e p)}.&(4.4)}$$

 We need
\proclaim L EMMA 4.2.  Suppose that $f\in {\cal M}  (\X )$, the increments
$d_j$'s are
$F_{j-1}$-conditionally symmetrically distributed, $j=1,2,\ldots,$
$\lambda
>0$, $\delta_1>0$, $\delta_2>0$, $\beta-1-\delta_2>0$.  Then
$$\P (f^*>\beta\lambda, w^* \le\lambda) \le \varepsilon
\P (f^*>\lambda),$$
where 
$$w^* = \left({d^*\over\delta_2}\right) \vee \left({D\,s_2\over
\delta_1}\right),$$
 $$\varepsilon = 2 \left({\e \over N}
{\delta^2_1\over
\delta^2_2}\right)^N, \quad N = {\beta-1-\delta_2\over \delta_2}.
\leqno
(4.5)$$

P{\eightrm ROOF.} Put $\bar d_j = d_j I\{\|d_j\| \le\delta_2
\lambda\}$,
$\bar f_j = \sum\limits^j_{i=0} \bar d_i$, $h_j = \bar f_{(j\wedge \tau \wedge \nu)\vee \mu} -\bar f_\mu,
$
 where
$\mu = \inf\{j\colon
\ \|\bar f_j\| > \lambda\}$, $\nu = \inf\{j\colon \ \|\bar f_j\| >
\beta\lambda\}$, $ \tau = \inf \{j\colon \ \bar s_{j+1} >
\lambda\}$,
$\bar s_{j+1} = \sqrt{\sum^{j+1}_{i=1} \E _{i-1} \|\bar
d_i\|^2}$, $j=0,1,\ldots\ $.
Then $\E _{j-1} \bar d_j = 0$, $h_j - h_{j-1} = \bar d_j I\{\mu<j\le
\tau\wedge \nu\}$, $j=1,2,\ldots$, and so, $(h_j)$ is a martingale in
$\X $
conditionally on $F_\mu$ [the $\sigma$-field consisting of all
$\Omega_0\in
F$ such that $\Omega_0 \cap \{\mu =j\} \in F_j$ for all $j$].
Besides,
$$\eqalign{\P (f^*>\beta\lambda,	w^*\le \lambda) &= \P (\bar f^*
>\beta\lambda, w^*\le \lambda)
\le \P (h^* > (\beta-1-\delta_2)\lambda)\cr
& = \E\P (h^* > (\beta-1-\delta_2)\lambda\mid F_\mu) I\{\mu<\infty\}\cr
&\le \varepsilon \P (\mu<\infty) = \varepsilon \P (f^*>\lambda);}$$
here, we have applied (3.5) with $r = (\beta-1-\delta_2)\lambda$,
$a=\delta_2\lambda$, $b=\delta_1\lambda$. \sqr

P{\eightrm ROOF OF} T{\eightrm  HEOREM} 4.1.  Argument in
Hitczenko
(1990) shows that we need to consider only  the following two
cases: 1) \ the increments $d_j$'s are $F_{j-1}$-conditionally
symmetrically
distributed, and  2) $f\in {\cal M}  _{ind}(\X )$. But via the standard
symmetrization formula $\widetilde X  = X -X '$, where $X ,X '$ are
independent
copies, one can easily reduce case 2) to 1). Thus, one can use
Lemma~4.2.

Applying now Lemma 7.1 in Burkholder (1973) with $\Phi(\lambda) =
\lambda^p$, $\gamma=\beta^p$, $\delta=\eta=1$, one has
$$\|f^*\|_p \le 2\beta \left({\|  d^*\|  _p\over \delta_2} +
{D\|  s_2\|  _p
\over\delta_1}\right) \quad \hbox{if}\quad \beta^p \varepsilon \le
1/2,
\leqno (4. 6)$$
where $\vp$ is given by (4.5). Choose now, for any $c\in [1,p]$,
$$\beta = 1+\e ^{-p/c}  + {1\over c}, \quad \delta_2= {1\over 10c},
\quad
\delta_1 = {1\over 10\sqrt c\, \e ^{p/c}}.$$
Then $\beta<3$,
$$\eqalign{N\left( = {\beta-1-\delta_2\over \delta_2}\right) =
9+10c\e ^{-p/c} >9,\cr
\beta^p = [1+(N+1)\delta_2]^p < \e ^{2pN\delta_2},\cr
(\beta^p\vp)^{1/N} < 2^{1/9} \e ^{2p\delta_2} {\e\over 10 c\e ^{-p/c}}
{c\over
\e ^{2p/c}} < {1\over 2},}$$
so that (4.6) implies
$$\|  f^*\|  _p \le 60(c\|  d^*\|  _p + \sqrt c\, \e ^{p/c}
D\|  s_2\|  _p);$$
thus, (4.1) is proved.

Let $c_p$ stand for the unique solution to the equation $\sqrt{c_p} =
\e ^{p/c_p}$. Then, $c_p \sim 2p/\ln p$ as $p\to \infty$. Hence, putting
$c=p$,
$c=c_p$, in (4.1), one comes to (4.2), (4.3), resp.

The function $g(c) := \sqrt c\, \e ^{p/c}$ decreases on $[1,p]$, $g(1)
=
\e ^p$, $g(p) = \e\sqrt p < \e p$. Hence, for each $\alpha \in [1, p/\ln (\e p)]$,
there exists $z_\alpha \in [1,p]$, the unique solution to the
equation
$g(z_\alpha) = \e ^{p/\alpha}$; besides, $g(2\alpha) =
\e ^{p/\alpha}\sqrt{2\alpha}\, \e ^{-p/(2\alpha)} \le \e ^{p/\alpha}
\sqrt{2\alpha_p} \e ^{-p/(2\alpha _p)} < g(z_\alpha)$, where $\alpha_p
:=
p/\ln (\e p)$. Thus, $z_\alpha <2\alpha$. Now we see that (4.1) with $c
=z_\alpha$ yields (4.4). \sqr

For $\X =\Re
%{\Bbb R}
$, inequalities (4.3) and (4.2) were proved in
Hitczenko
(1990) and Hitczenko (1991), resp.

A spectrum of bounds on moments of martingales in Hilbert spaces with
bounded second conditional moments was found in Pinelis (1980). It is
essentially equivalent to (4.1) at least in the case of independent
increments $d_j$'s (see Proposition~9.2 below) but has a much more
cumbersome expression.

The infimum in $c$ of the R.H.S. of (4.1), evaluated in section 6 below,
turns out to be  an
upper bound on $\|  f^*\|  _p$ optimal in terms of $\|  d^*\|  _p$,
$\|  s_2\|  _p$, the optimum choice of $c$ depending, obviously, on
$\|  d^*\|  _p/\|  s_2\|  _p$. Besides, for each $c\in [1,p]$, the
``individual''
bound $c\|  d^*\|  _p + \sqrt c\, \e ^{p/c}\|  s_2\|  _p$ is optimal for a
certain
corresponding value of $\|  d^*\|  _p/\|  s_2\|  _p$; in particular, all
the
bounds in (4.2)-(4.4) are optimal. The issue of optimality is treated
rigorously in section 6.

\proclaim T HEOREM 4.3. If $\X $ is a $(2,D)$-smooth separable Banach space,
$f\in
{\cal M}  (\X )$, and the increments $d_j$'s are $F_{j-1}$-conditionally
symmetrically
distributed, then
$$\|  f^*\|  _p \< \sqrt p\, D\|  S_2\|  _p,\qquad p\ge 1. \leqno $$

P{\eightrm ROOF.} Consider $f$ as a martingale relative to the sequence
of
$\sigma$-fields $(G_j)$ defined in the proof of Theorem~3.6. Then
$s_2=S_2$.
Reasoning as in the proof of Lemma~4.2 but using Theorem~3.6 instead
of
Theorem~3.4, we see that
$$\P (f^*>\beta \lambda, D\,S_2\le\delta\lambda) \le \vp \P (f^*>\lambda), \leqno (4.7)$$
$\lambda>0, \delta>0,\, \beta-1-\delta>0$, where $\vp =
\exp[-\delta^{-2}(\beta-1-\delta)^2/2]$. It remains to choose, say,
$\beta=2$, $\delta=(0.1) p^{-1/2}$ and apply Lemma~7.1 of Burkholder
(1973). \sqr

\n R{\eightrm EMARK.} For conditionally symmetric martingales in Hilbert spaces,  the exact constant $A_p$ in the inequality 
$$\|  f_n\|  _p \le A_p \|  S_2\|  _p\leqno (4.8)$$
 was found [see
Burkholder
(1991, chapter~11), Davis (1976), Wang (1991)]. For any real
martingale
with independent symmetrically distributed increments, it follows
from the
result of Whittle (1960) [proved correctly only for $p=2$ and $p\ge 3$] and of Haagerup (1982) that (4.8) takes place with  the exact, in this "independent increments" case, constant $A_p =
\|  \xi\|  _p$, where $\xi \sim N(0,1)$,
 $p\ge 2$.
  Since $\|  \xi\|  _p \asymp
\sqrt p$,
the bound in Theorem~4.3 is optimal (to $\asymp$). \sqr 

\n R{\eightrm EMARK.} Bounds  given in  Theorems 3.1,    4.1, 4.3 are only possible in the {\it 2-smooth\/}  Banach spaces, even if we need a bound like those in Theorem 4.1 for just one particular $p$. Indeed, all the mentioned here results imply $\|f^*\|_p \le C  D \|  S_2\|  _p $ for at least one particular $p\ge 2$, some $C>0$  and all Walsch-Paley martingales in $\X$ since, for those martingales, $s_2=S_2\ge d^*$. Thus, one has (4.7) with $\vp = (C D \delta/(\beta -1-\delta))^p$, and so, $\| f^*\|_2 \le C_1  D \|  S_2\|  _2 $ for some $C_1>0$. It remains now to recall the characterization of the 2-smooth Banach spaces given by Pisier (1975). \sqr

\vskip .8cm

\n {\bf 5. Applications:\ bounds on central moments of the norm of
the sum
of independent random vectors in arbitrary Banach spaces.}

\proclaim T HEOREM 5.1. If $f\in {\cal M}  _{ind}(\X )$, $(\X, \|\cdot\|) $ is any separable
Banach
space, $p\ge 2$, $1\le c\le p$, and $x$ is any non-random vector in $\X $,
then for all $n=1,2, \ldots $,
$$\|  \, \|f_n+x\| - \E \|f_n+x\|\, \|  _p \< c\|  d^*\|  _p + \sqrt c\,
\e ^{p/c}
\|S_2\|_2;\leqno (5.1)$$
in particular,
$$\leqalignno{\|  \, \|f_n +x\| - \E \|f_n+x\|\,\|  _p &\< p \|  d^*\|  _p +
\sqrt
p\, \|S_2\|_2,&(5.2)\cr
\|  \, \|f_n+x\| - \E \|f_n+x\|\, \|  _p &\< {p\over \ln p} (\|  d^*\|  _p
+\|S_2\|_2),&(5.3)\cr
\|  \, \|f_n+x\| - \E \|f_n+x\|\, \|  _p &\< \alpha \|  d^*\| _p + \e ^{p/\alpha}
\|S_2\|_2,&(5.4)}$$
$1\le \alpha \le p/\ln (\e p)$.

P{\eightrm ROOF.} We use the following modification of the method of
Yurinskii (1974), being stated here as

\proclaim L EMMA 5.2. For $f\in {\cal M}  _{ind}(\X )$, let $F_j$ be the $\sigma$-field generated by $f_0,f_1, \ldots , f_j$; put $\zeta_j =
\E  _j\|f_n+x\| -
\E \|f_n+x\|$, $\xi_{j+1} = \zeta_{j+1} - \zeta_j$, $j=0,1,\ldots, n$.
Then
$(\zeta_j)$ is a martingale,
$$|\xi_j| \le \|d_j\| + \E \|d_j\|,\quad \E  _{j-1} \xi^2_j \le
\|d_j\|^2,$$
 $j=1,\ldots, n$, and
$$\|f_n+x\| - \E \|f_n+x\| = \zeta_n.$$

P{\eightrm ROOF}  can be found in each of the papers of Pinelis (1981, 1986),
Pinelis
and Sakhanenko (1985). Since it is very short, we repeat it here for
the
reader's convenience. Put $\eta_j = \E  _j(\|f_n+x\| - \|f_n+x-d_j\|)$.
Then
$|\eta_j| \le\|d_j\|$ .  Besides, $\xi_j = \eta_j-\E  _{j-1}\eta_j$,
$j=1,\ldots, n$, in view of the  independence of $d_j$'s; hence,
$$\eqalign{|\xi_j|&\le |\eta_j| + |\E  _{j-1}\eta_j| \le \|d_j\| +
\E \|d_j\|,\cr
\E  _{j-1}\xi^2_j &= \E  _{j-1}\eta^2_j - (\E  _{j-1}\eta_j)^2 \le \E  _{j-1}
\eta^2_j
\le \E \|d_j\|^2,\quad j=1,\ldots, n. \sqr
}$$ 

To complete now the proof of Theorem~5.1, apply Theorem 4.1 to the
martingale $(0,\zeta_1,\zeta_2,\ldots,
\zeta_n,\zeta_n,\zeta_n,\ldots)$,
$\zeta_j$'s defined in Lemma~5.2. Using the lemma and the
inequalities of
Minkowsky and H\"older, we see that
$$\eqalign{\sum^n_{j=1} \E  _{j-1} \xi^2_j &\le \sum^n_{j=1} \E \|d_j\|^2
=
\|S_2\|_2^2,\cr
\|  \max_{j\le n}|\xi_j|\, \|  _p &\le \|  \max_j \, (\|d_j\| +
\E \|d_j\|)\|  _p\cr
&\le \|  d^*\|  _p + \max_j \E \|d_j\|\cr
&\le 2\|  d^*\|  _p.}$$
Finally, note that the constant of smoothness $D(\Re )=1$. Now the  statement of Theorem~5.1 is
evident. \sqr

Hoffmann-J\o rgensen (1974) found the following extension of the
Rosenthal's
inequality for $f\in {\cal M}  _{ind}(\X )$, $p\ge 1$:
$$\| \; \|  f_n\| \; \| _p \le c(p)( \|  d^*\|  _p + \| \; \|  f_n\| \; \| _1)\leqno (5.5)$$
with $c(p)$ depending only on $p$ [it can be seen  that the best choice of
parameters in the method of  Hoffmann-J\o rgensen gives (5.5) with
$c(p) \asymp
p$].

For $f\in {\cal M}  _{ind}(\X )$, Pinelis (1978) proved that
$$\| \; \|  f_n\| \; \|  _p \< c_1(p) \| \; \|  f_n\| \; \|  _1 + c_2(p) \|S_p\|_p + c_3(p)\|S_2\|_2,$$
which is also a generalization of the Rosenthal's (1970) inequality;
the
method can actually yield $c_1(p) =1$, $c_2(p) = p$, $c_3(p) = \sqrt
p$.

An inequality, implying (5.1), was obtained in Pinelis (1980) [see
Proposition~9.2 below].

De Acosta (1981) proved a version of (5.3) but with $x=0$ and with an
implicit factor $c(p)$ instead of $p/\ln p$. Since he used 
inequality (21.1) 
of Burkholder (1973), Proposition~9.1 below shows that in fact 
his
method leads to $c(p) = p$.

Using isoperimetric technique, Talagrand (1989) proved the following
version
of (5.5):
$$\| \; \|  f_n\| \; \|  _p \< {p\over \ln (2p)} (\|  d^*\|  _p + \| \; \|  f_n\| \; \|  _1) \leqno
(5.6)$$
for $f\in {\cal M}  _{ind}(\X )$, $p\ge 1$, which was also proved in Kwapie\'n
and
Szulga (1991) by a different method. The Talagrand's inequality may
be
compared with (5.1) and, in particular, with (5.3). If $\|S_2\|_2$ is much
greater than $\|  d^*\|  _p + \| \; \|  f_n\| \; \|  _1$, then (5.1) may lose to
(5.6) in
certain cases. If however $\X $ is e.g. of cotype~2, then (5.3) is often no
worse
then (5.6); moreover, say, in the typical case of the increments with
the
same or almost the same distribution, (5.2) is significantly better.
On the
other hand, if there are heavy distribution tails, i.e., if
$\|  d^*\|  _p$
is much greater than both $\|S_2\|_2$ and $\| \; \|  f_n\| \; \|  _1$, then (5.4) with
$\alpha
\asymp 1$ does better than (5.6). Other advantages of bounds
(5.1)-(5.4)
are that they are applicable to the sums of non-zero-mean random
vectors
(owing to the presence of $x$) and better reflect the concentration
phenomenon of the distribution of the sum of independent random
vectors.

Modifications of the method of Yurinskii (1974) for $f\in {\cal M}  _{ind}(\X )$, allowing reduction of the problem of
upper
bounds on the L.H.S. of (5.1) for any separable Banach space $\X $ to
that of
upper bounds on $\|  f_n\|  _p$ for $f\in {\cal M}  _{ind}(\Re
%{\Bbb R}
)$ were
proposed
in Berger (1991) [for $x=0$] and in Pinelis (1992a) [with the best constant, for any $x\in \X$];
actually,
instead of the power moment function $u\mapsto |u|^p$, one can use
any convex function there.

A straightforward application of Theorem~4.1 yields the following
bounds in
the case of sums of independent zero-mean random variables.

\proclaim T HEOREM 5.3. Suppose that $f\in {\cal M}  _{ind}(\X )$, $\X $ is a
Hilbert
space, $p\ge 2$, $1\le c\le p$. Then
$$\| \; \|  f_n\| \; \|  _p \< c\| S_p \| _p + \sqrt c\, \e ^{p/c} \|S_2\|_2;\leqno (5.7)$$
in particular,
$$\leqalignno{\| \; \|  f_n\| \; \|  _p &\< p\| S_p \| _p + \sqrt p\, \|S_2\|_2,&(5.8)\cr
\| \; \|  f_n\| \; \|  _p &\< {p\over \ln p} (\| S_p \| _p+\|S_2\|_2),&(5.9)\cr
\| \; \|  f_n\| \; \|  _p &\< \alpha \| S_p \| _p + \e ^{p/\alpha} \|S_2\|_2, \qquad 1\le \alpha  \le
{p\over \ln (\e p)}.&(5.10)}$$

P{\eightrm ROOF},  in view of Theorem 4.1 and the equality $D(\X)=1$ for any Hilbert space $\X$, consists in the
following
trivial remarks:\ $\| \;  \|f_n\| \; \| _p \le \|  f^*\|  _p$,
$\|  d^*\|  _p \le \| S_p \| _p, $ 
and, for $f\in {\cal M}  _{ind}(\X )$, $s_2=S_2$.\sqr

As it was said in section 4, the results of Pinelis (1980) imply
(5.1),
(5.7); it was also explained in that paper how to elicit bounds like
(5.8)
and (5.10) (for $\alpha=4$). Nevertheless, it is not obvious how to
deduce
a general inequality like (5.7) from the spectrum of bounds in
Pinelis
(1980) (again, we refer to Proposition~9.2). From this point of
view, even
in the classical case of sums of independent real-valued zero-mean
random
variables, (5.7) is apparently new.

An inequality similar to, say, (5.9), but with $2^p$ instead of
$p/\ln p$,
was probably first found by Rosenthal (1970) for $f\in {\cal M}  _{ind}(\Re
%{\BbbR}
)$,
who also obtained a lower bound, which differs from the upper one by at most 
a
factor depending only on $p$.

Ros\'en (1970) proved a result for $f\in {\cal M}  _{ind}(\Re
%{\Bbb R}
)$, which
implies the upper Rosenthal's bound for $p=2,4,6,\ldots$ [this
implication
was unnoticed; to demonstrate it, one can put, in the notation of
Ros\'en,
$\lambda_\nu(p) = (\E X ^{2p}_\nu/\E X ^2_\nu)^{1/(2p-2)}$, $
\rho_\nu(p) =
\E X ^2_\nu/(\lambda_\nu(p)^2)$].
Moreover, using some ideas of Dharmadhikari and Jogdeo (1969), it
is
possible to deduce the upper bound by Rosenthal for all real $p\ge 2$ from
the
Ros\'en's result.

A method like just described was  used in the student diploma work of
Pinelis (1974) to prove an upper bound of the Rosenthal type for
$f\in
{\cal M}  _{ind}(\Re
%{\Bbb R}
)$ [via the  Marcienkiewicz-Zygmund (1937)
inequality, the
lower Rosenthal bound was also obtained therein]. While the constants
in
Pinelis (1974) were implicit, the method could yield (5.9);
regrettably,
the results of Rosenthal (1970) and Burkholder (1973) had long
remained
unknown to the author, and so, the problem of the constants was not
among
the ones considered then by him as most urgent.

Inequality (5.9) for  $f\in {\cal M}  _{ind}(\Re
%{\Bbb R}
)$ and for sums of
exchangeable
random variables, with the proof that $p/\ln p$ is optimal in (5.9),
was
first given by Johnson, Schechtman and Zinn (1985).

A.I. Sakhanenko, a referee of the mentioned diploma work, upon having
reviewed it, suggested another approach, giving in effect (5.8),
again for
$f\in {\cal M}  _{ind}(\Re
%{\Bbb R}
)$ [see Nagaev and Pinelis (1977)].

An inequality, similar to (5.10) with $\alpha=4$, was found by
Sazonov
(1974).

All bounds in Theorems~4.1, 5.1, 5.3 are optimum. We shall prove the
optimality in the next section, using results and ideas from Pinelis
and
Utev (1984), where, in particular, for any $p\ge 2$,
$$\sup\{\|  f_n\|  _p\colon \  f\in {\cal M}  _{ind}(\Re
%{\Bbb R}
),\ \| S_2 \| _2\ 
\hbox{and}\ \| S_p \| _p \ \hbox{fixed}\}$$
was computed up to $\asymp$, which, for instance,
implies all inequalities (5.7)-(5.10) for $f\in {\cal M}  _{ind}(\Re
%{\Bbb R}
)$;
also,
it was noted in Pinelis and Utev (1984) that bounds like (5.8) and
(5.10)
[for $\alpha \asymp 1$] represent in a certain sense the two extreme
bounds
in the spectrum of all optimal bounds on moments, the optimum value
of
a ``spectrum parameter'' depending on $\| S_p \| _p/\| S_2 \| _2$. 

\vskip .8cm

\n {\bf 6. Optimality of the bounds on moments.}
\ \ Let us consider the following upper bounds for any $a_2>0$,
$a_p>0$, $p\ge 2$:
$$\eqalign{B_p := B_p (a_p,a_2)
:= \sup\{&\|  f^*\|  _p\colon \ f\in {\cal M}  (\X ),\ \X 
\hbox{ is $(2,D)$-smooth,}\cr
& \|  d^*\|  _p = a_p,\ D\|  s_2\|  _p = a_2\},\cr
B^\X _{p,{ind}} := B^\X _{p,ind}(a_p,a_2):= \sup\{&\|  \, \|f_n+x\| - \E\|f_n+x\|\, \|  _p \colon \ f\in
{\cal M}  _{ind}(\X ),\cr
& x\in \X , \ \|  d^*\|  _p = a_p,\ \|S_2\|_2=a_2, \ n=1,2, \ldots\},\cr
B_{p,{ind},0} := B_{p, {ind}, 0}(a_p,a_2):= \sup\{&\|  f_n\|  _p \colon \ f\in {\cal M}  _{ind}(\Re
%{\Bbb R}
),\cr
& \|  d^*\|  _p = a_p,\ \|S_2\|_2=a_2, \ n=1,2, \ldots\},}$$
and their analogues $B_{p,S} = B_{p,S}(a_p,a_2)$,
$B^\X _{p, {ind}, S} = B^\X _{p, {ind}, S}(a_p,a_2)$, $B_{p, {ind},
0, S} \break =
B_{p, {ind}, 0, S} (a_p,a_2)$
obtained by replacing the equality $\|  d^*\|  _p = a_p$ in the above
definitions by $\|S_p\|_p=a_p$.

We shall show that all the introduced bounds are $\asymp$-equivalent to each
of the
following:
$$\eqalign{B^*_p &:= B^*_p (a_p,a_2) := a_p + \sqrt p\, a_2 +
{pa_p\over \ln
(2+{a_p\over a_2} \sqrt p)},\cr
\hat  B_p &:= \hat  B_p (a_p,a_2) := \min\{ca_p + \sqrt c\,
\e ^{p/c}
a_2\colon \ 1\le c\le p\},}$$
$$\check B_p := \check B_p (a_p,a_2) := \max\{(p\alpha
+1)^{1-\alpha/2}
a^{1-\alpha}_p a^\alpha_2\colon \ 0 \le \alpha \le 1\}.$$

Theorem 6.1 below principally means that for any pair
$(a_p,a_2)$ of
the values of the characteristics  used in the bounds (4.1), (5.1)
and
(5.7) there exists a value of the ``spectrum parameter'' $c$ providing an
optimum bound. Roughly speaking, it means that spectra (4.1), (5.1)
and
(5.7) are rich enough, so Theorem 6.1 may be called the ``spectrum
completeness
theorem''. It also means that it is not essential in this context values of
which of
the two pairs are fixed:\  $\|  d^*\|  _p$ and $\|  s_2\|  _p$,
or
 $\|S_p\|_p$ and $\|  s_2\|  _p$.

\proclaim T HEOREM 6.1. For all $p>2$, $a_2>0$, $a_p>0$ and for all
separable Banach spaces $(\X , \|\cdot\|)$,
$$\eqalign{B_p \asymp B^\X _{p,{ind}} &\asymp B_{p, {ind},0}  \asymp
B_{p,S}
\asymp B^\X _{p,{ind},S} \asymp B_{p,{ind},0,S}\cr
&\asymp \check B_p \asymp \hat  B_p \asymp B^*_p.}$$

Proof is comparatively long and will be given later in this section.

The ``spectrum'' $ca_p + \sqrt c\, \e ^{p/c}a_2$, $1\le c\le p$, turns
out to
be not only ``complete'' but also ``minimal'' in the sense that for
each $c
\in [1,p]$, there exist $a_p>0$, $a_2>0$ such that the ``individual''
bound
$ca_p + \sqrt c\, \e ^{p/c}a_2$ is the best possible. Let us now give
the
rigorous statement.

\proclaim T HEOREM 6.2. For any $c\in[1,p]$, $p>2$,
$$\sup\left\{ {B_p (a_p,a_2)\over ca_p + \sqrt c\, \e ^{p/c}a_2}\colon \
a_p
>0, a_2>0\right\} \asymp 1;$$
here any of the other 5 bounds $B^\X _{p,{ind}}, B_{p,{ind},0}$,
$B_{p,S}$,
$B^\X _{p,{ind},S}$, $B_{p,{ind},0,S}$ may be used in place of $B_p$. In
particular, all bounds (4.2)-(4.4), (5.2)-(5.4), (5.8)-(5.10) are
optimal.

Proof will be given after that of Theorem 6.1.

The following proposition might seem analogous to Theorem~6.2, but it
is
less important since, for an ``individual'' value of $\alpha$,
$(p\alpha
+1)^{1-\alpha/2} a^{1-\alpha}_pa^\alpha_2$ does not represent an
upper
bound on the moments -- actually, by Theorem~6.1, $B_p (a_p,a_2) \>
(p\alpha+1)^{1-\alpha/2}a^{1-\alpha}_pa^\alpha_2$.

\proclaim P ROPOSITION 6.3. For any $\alpha \in [0,1], p>2$,
$$\inf\left\{{B_p (a_p,a_2)\over
(p\alpha+1)^{1-\alpha/2}a^{1-\alpha}_p
a^\alpha_2}\colon \ a_p>0, a_2>0\right\}\asymp 1;$$
here, any of the other 5 bounds $B^\X _{p,{ind}}$, $B_{p,{ind},0}$, 
$B_{p,S}$,
$B^\X _{p,{ind},S}$, $B_{p,{ind},0,S}$ may be used in place of $B_p$.

Proof will be given after that of Theorem 6.2.

\n R{\eightrm EMARK.} It is easy to see that  Doob's inequality
$$\|  f^*\|  _p \le {p \over {p-1}}\, \sup _n \| \; \| f_n\| \; \| _p\,  , \quad p>1, f\in {\cal M}  (\X ),$$
(see also (1.4) in Burkholder (1973)) remains true for all separable
Banach
spaces $\X $. Therefore, one could replace $\|  f^*\|  _p$ in the
definition of
$B_p$ by $\sup _n  \| \; \| f_n\| \; \| _p$, and statements \break 6.1--6.3 would hold. This
remark can
be also deduced from the proof of these statements. \sqr

Proof of Theorem 6.1 is based on  lemmas 6.4 -- 6.9 below.

\proclaim L EMMA 6.4. If $p>2$, $m$ is a positive integer, $m\le p/2$,
$a_p>0$, $a_2>0$, then
$$B_{p,{ind},0} (a_p,a_2) \> a_p+B_{2m,{ind}, {sym}, S}((a^{2m-2}_p
a^{p-2m}_2)^{1\over p-2} , a_2),$$
where $B_{2m,{ind}, {sym}, S}$ may be defined by (6.6), (6.7) below.

P{\eightrm ROOF.} Put $g_n(t) = t^{-p/2} [1-(1-t/n)^n]$, $g_\infty(t) =
t^{-p/2}(1-\e ^{-t})$,  $t>0$, $n>0$. Then $g_n(t)$
decreases
in $n$ to $g_\infty(t)$, $g_n(t)$ decreases in $t\in (0,n]$ from
$g_n(0+) =
\infty$ to $g_n(n) = n^{-p/2}$, $g_\infty(t)$ decreases in $t\in
(0,\infty)$ from $g(0+) = \infty$ to $g(\infty) = 0$. Hence, for any
$n\in(a^2_2/a^2_p, \infty]$, there exists the unique solution $t_n\in
(0,n)$
to the equation $g_n(t_n) = (a_p/a_2)^p$, and $t_n\downarrow t_\infty$ as
$n\uparrow
\infty$.

For any $p>0, a_p>0, a_2>0$, and any integer $n>a^2_2/a^2_p$,
consider
$f\in {\cal M}  _{ind}(\Re
%{\Bbb R}
)$ such that
$$\P (d_j=\pm u_n) = q_n/2,\quad \P (d_j=0) = 1 -q_n \quad\hbox{for}\quad
j=1,\ldots, n,$$
$$\P (d_j=0) = 1\quad \hbox{for}\quad j=n+1, n+2,\ldots,$$
where
$$q_n := {t_n\over n},\quad u_n := {a_2\over \sqrt{t_n}}.$$
Then, for the so constructed $f$,
$$\|  d^*\|  _p = a_p,\quad \|S_2\|_2 = a_2,\leqno(6.1)$$
$$\|S_{2m}\|_{2m} = t^{(1-m)/(2m)}_n a_2\to t^{1/(2m)}_\infty
u_\infty,\leqno(6.2)$$
$$\E\e ^{i\lambda f_n} = [1+q_n(\cos \lambda u_n-1)]^n \to \E\e ^{i\lambda
Z},
\quad \lambda\in \Re
%{\Bbb R}
,$$
as $n\to \infty$, where
$$u_\infty = \lim_{n\to \infty} u_n = {a_2\over
\sqrt{t_\infty}},\leqno
(6.3)$$
$Z$ is a (symmetrized Poisson) random variable with the
characteristic
function
$$\E\e ^{i\lambda Z} = \exp[(\cos \lambda u_\infty  - 1)t_\infty].\leqno
(6.4)$$
Hence, by the analogue of the Fatou lemma for the convergence in
distribution [see, e.g., Theorem~5.3 in Billingsley (1968)],
$$\liminf_{n\to \infty} \|  f_n\|  _{2m} \ge \|  Z\|  _{2m}.\leqno
(6.5)$$
Consider now, for $b_{2m}>0, b_2>0$,
$$\eqalign{B_{2m,{ind}, {sym}, S}(b_{2m}, b_2) :=
&\sup\{\|  f_n\|  _{2m}\colon \ f\in {\cal M}  _{ind}(\Re
%{\Bbb R}
),\  \|S_{2m}\|_{2m}=b_{2m},\ 
\|S_2\|_2=b_2,\cr
&\qquad d_j\hbox{'s}\hbox{ symmetrically distributed }\}.}$$
The theorem in Pinelis and Utev (1984) states that
$$B_{2m, {ind}, {sym}, S}(b_{2m},b_2) =
\left(\sum^m_{j=1}\Gamma_{j,m}
(b^{m-j}_{2m}b^{j-1}_2)^{2m\over{m-1}}\right)^{1\over 2m},\leqno (6.6)$$
where
$$\Gamma_{j,m} = (2m)! \sum^j_{r=1} \sum \prod^r_{k=1}  {1 \over {
j_k!((2m_k)!)^{j_k}  }  };\leqno (6.7)$$
the inner summation is over all positive integers $m_1>\cdots > m_r$
and
$j_1,\ldots, j_r$ such that $m_1j_1+\cdots + m_rj_r= m$, $j_1+\cdots
+j_r=j$. It follows from its proof that
$$B_{2m,{ind}, {sym,}S}(t^{1\over 2m}_\infty u_\infty, t^{1\over
2}_\infty
u_\infty) = \|  Z\|  _{2m},\leqno (6.8)$$
where $Z$ is defined by (6.4) [note the misprint in formula (7) of
Pinelis
and Utev (1984): there must be $A^{j_k}_{2m_k,n}$ instead of
$A_{2m_k,n}$].
Observe that for the $f$ constructed above in this proof,
$$\|S_{2m}\|_{2m} = (\|S_p\|_p^{2m-2} \|S_2\|_2^{p-2m})^{1\over{p-2}}.$$
Now, (6.2), (6.1) and the inequality $\|S_p\|_p\ge \|  d^*\|  _p$ imply
$$t^{1\over 2m}_\infty u_\infty\ge (a^{2m-2}_p
a^{p-2m}_p)^{1\over{p-2}}.\leqno
(6.9)$$
But (6.6), (6.7) show that $B_{2m, {ind}, {sym}, S}$ is an increasing
function in either of its arguments. Thus, using (6.5), (6.8), (6.9)
and
(6.3), we deduce
$$\liminf_{n\to \infty} \|  f_n\|  _{2m} \ge B_{2m, {ind}, {sym}, S}
\( (a^{2m-2}_p a^{p-2m}_2)^{1\over p-2}, a_2 \) .$$
The lower Rosenthal's bound and   H\"older's inequality give
$$\|  f_n\|  _p \> \|S_p\|_p \vee \|  f_n\|  _{2m} \ge a_p \vee \|  f_n\|  _{2m}.$$
Now the statement of the lemma is obvious. \sqr

\proclaim L EMMA 6.5. Under the conditions of Lemma 6.4,
$$B_{p, {ind}, 0,S}(a_1,a_2) \> a_p +B_{2m, {ind}, {sym},
S}((a^{2m-2}_p
a^{p-2m}_2)^{1\over p-2}, a_2).$$

P{\eightrm ROOF}  is similar to that of Lemma 6.4 but easier. In fact, Lemma~6.5
was
proved in Pinelis and Utev (1984). \sqr

\proclaim L EMMA 6.6. For all $p\ge 2$, $a_p>0$, $a_2>0$, one has
$\hat 
B_p (a_p,a_2) \< \check B_p (a_p,a_2)$.

P{\eightrm ROOF.} Consider the functions
$$\leqalignno{g_p (\alpha) &:= (p\alpha+1)^{1-\alpha/2}q^\alpha, \qquad q :=
a_2/a_p, \qquad 0 \le \alpha \le 1,&(6.10)\cr
q_p (\alpha) &:= \sqrt{p\alpha +1} \,  \exp\left\{{1\over 2} {p\alpha
-2p\over  p\alpha+1}\right\},\qquad 0 \le \alpha \le 1.&(6.11)}$$
Then, by the definition of $\check B_p$,
$$\check B_p (a_p,a_2) = \max\{a_pg_p (\alpha)\colon \ 0 \le \alpha \le 1\}.\leqno (6.12)$$
Besides, $q_p$ is continuous and increasing and maps [0,1] onto the segment
$$I_p := \left[\e ^{-p}, \sqrt{p+1} \, \exp \left\{ -{1\over 2} 
{p\over{p+1}}\right\}\right]. \leqno (6.13)$$

An essential relation between $g_p$ and $q_p$ is
$$g'_p (\alpha) = g_p (\alpha) \ln {q\over q_p (\alpha)}\, .\leqno (6.14)$$

If we have $q\in I_p$, take $\alpha_p = q^{-1}_p (q)$; then $\alpha_p
\in
[0,1]$, $q=q_p (\alpha_p)$, and, in view of (6.10), (6.11),
$$\leqalignno{(p\alpha_p+1)a_p &= a_pg_p (\alpha_p)
\exp\left\{-{\alpha_p\over 2} {p\alpha_p-2p\over p\alpha_p
+1}\right\} \asymp a_pg_p (\alpha_p),&(6.15)\cr
\sqrt{p\alpha_p+1}\, \e ^{p/(p\alpha_p+1)} a_2 &= a_pg_p (\alpha_p)
\exp\left\{{p\alpha_p (3-\alpha_p)\over 2(p\alpha_p+1)}\right\} \asymp a_pg_p (\alpha_p).}$$
This implies, with $c_p := p\alpha_p+1$, that
$$c_pa_p + \sqrt{c_p}\, \e ^{p/c_p} a_2 \asymp a_pg_p (\alpha_p)\leqno
(6.16)$$
(if $q\in I_p$). Obviously, $c_p\in [1,p+1]$. Hence, in view of
(6.12), if
$q\in I_p$ and $c_p\le p$, the lemma is true; if however  $c_p\in
(p,p+1]$,
then
$$\eqalign{pa_p + \sqrt p\, \e ^{p/p} a_2 &\le c_pa_p + \sqrt{c_p}\,
\e ^{p/c_p} \e a_2\asymp a_pg_p (\alpha_p)}$$
by (6.16). Thus, the lemma is true whenever $q\in I_p$.

Consider now the two cases when $q\notin I_p$.

First, suppose that
$$q < \min I_p = \e ^{-p}.\leqno (6.17)$$
Then $q < q_p (\alpha)$ for all $\alpha \in [0,1]$, and (6.14) implies
$g'_p
<0$ on $[0,1]$. Hence, $\check B_p (a_p,a_2) =a_pg_p (0)=a_p$; on the
other 
hand,
$\hat  B_p (a_p,a_2) \le a_p+\e ^pa_2=$ $(1+\e ^pq)a_p \le 2a_p$ by
(6.17).
Thus, the lemma is true in the case (6.17), too.

Finally, let
$$q>\max I_p (> \sqrt p /\e).\leqno (6.18)$$
Then $q>q_p (\alpha)$ for all $\alpha \in [0,1]$, and (6.11) implies
$\check
B_p (a_p,a_2) = a_p g_p (1)\ge \sqrt p\, a_2$, whereas $\hat  B_p
(a_p,a_2) \le
pa_p + \sqrt  p\, a_2=$ $(p/q +\sqrt p\, \e) a_2 \le$ $2\e\sqrt p\,
a_2$, in
view of (6.18). \sqr

\proclaim L EMMA 6.7. For all $p>2$, $a_p>0$, $a_2>0$, one has 
$B_{p,{ind},0}(a_p,a_2) \> \check B_p (a_p,a_2)$.

P{\eightrm ROOF.} By Lemma 6.4 and the estimate
$$\Gamma^{1/(2m)}_{j,m} \asymp j^{1-j/(2m)}\leqno (6.19)$$
[see (9) in Pinelis and Utev (1984)],
$$B_{p,{ind}, 0}(a_p,a_2) \> a_p+\widetilde B_p (a_p,a_2),\leqno
(6.20)$$
where
$$\widetilde B_p (a_p,a_2) = \max\{j^{1-j/2m} a^{1-\alpha_j}_p
a_2^{\alpha_j}\colon \ j=1,\ldots, m\},$$
$\alpha_j := 1-{2(m-j)\over p-2}$, $m$ is the integer part of $p/2$.
Since
$j\>p\alpha_j+1$ and $j/(2m) \le \alpha_j/2 + 1/p$, one has
$j^{1-j/2m} \>
(p\alpha_j+1)^{1-\alpha_j/2}$, $j=1,\ldots, m$. Therefore,
$$\widetilde B_p (a_p,a_2) \> \widetilde{\widetilde B}_p (a_p,a_2) :=
\max\{a_p g_p (\alpha_j)\colon \ j=1,\ldots, m\},\leqno (6.21)$$
where $g_p$ is defined by (6.10).

Consider first the case $q\in I_p$, where $q$ and $I_p$ are given by
(6.10)
and (6.13), resp. We can assume that $p>3$, because for $p\in [2,3]$,\quad  
$\check
B_p (a_p,a_2) \< a_p + a_2 \< B_{p,{ind}, 0}(a_p,a_2)$ in view of the
Rosenthal's lower bound. Then,
$$0 \le \alpha_1 < {2\over p-2} \< {1\over p},\qquad 0  < \alpha_{j+1}
-\alpha_j = {2\over p-2} \< {1\over p};$$
besides, $\alpha_m=1$. Note also that $q\in I_p$ implies $|\ln q|\<p$.
This,
along with (6.14), (6.11), leads to $|g'_p (\alpha)| \< pg_p (\alpha)$,
and
hence to $g_p (\alpha+\delta)\asymp g_p (\alpha)$ if $0 \le \alpha \le
\alpha+\delta \le 1$ and $\delta \< 1/p$. Thus,
$$\eqalign{\check B_p (a_p,a_2) = \max\{a_pg_p (\alpha)\colon \ 0 \le
\alpha
\le 1\}\asymp \max\{a_pg_p (\alpha_j)\colon \ j=1,\ldots, m\},}$$
and, in virtue of (6.20), (6.21), the lemma is true if $q\in I_p$.

Let us now consider the case (6.17). Then, as at the end of the proof Lemma 6.6, we see that  $\check B_p (a_p,q_2) =
a_pg_p (0)
=a_p$, and (6.20) implies the lemma.

Finally, suppose that (6.18) takes place. Then $\widetilde{\widetilde
B}_p (a_p,a_2) = a_pg_p (1) =  \check B_p (a_p,a_2)$, and now (6.20) and
(6.21) imply the lemma. \sqr

\proclaim L EMMA 6.8. For all $p>2, a_p>0, a_2>0,$ one has $ B_{p,{ind},0,S}
(a_p,a_2) \> \check B_p (a_p,a_2)$.

P{\eightrm ROOF}  differs from that of Lemma 6.7 only in that we have to refer to
Lemma~6.5 instead of Lemma~6.4. \sqr

\proclaim L EMMA 6.9. For all $p\ge 2, a_p>0, a_2>0$, one has
$B^*_p (a_p,a_2) \asymp \check B_p (a_p,a_2)$.

P{\eightrm ROOF.} Consider first the case $q\in I_p$, where $q$ and
$I_p$ are
defined in (6.10) and (6.13), resp. Then (6.15) implies
$$\check B_p (a_p,a_2) = a_pg_p (\alpha_p) \asymp c_pa_p,\leqno
(6.22)$$
where $g_p$ is given by (6.10), $\alpha_p = q^{-1}_p (q)$,
$c_p=p\alpha_p+1$. Putting
$z = (p+1/2)/ c_p$,
one has $z>1/2$. We can rewrite the equation $\alpha_p = q^{-1}_p (q)$ as
$\sqrt z\, \e ^z = {\sqrt{\e(p+{1/ 2})}\,a_p/a_2}$.
Since $\e ^z/2 < \sqrt z\, \e ^z < \e ^{2z}$ when $z>1/2$, we deduce
$z \asymp \ln (2+a_p\sqrt p/a_2)$,
$$c_p \asymp {p \over \ln (2+{a_p\over a_2}\sqrt p)}.\leqno
(6.23)$$
This implies $\check B_p \< B_p^*$ if $q\in I_p$.

In the cases (6.17) and (6.18), $\check B_p (a_p,a_2)$ equals
$a_pg_p (0) =
a_p$ and $a_pg_p (1) \asymp \sqrt p \, a_2$, resp. Taking into
account
(6.22), (6.23), (6.12), we see that in any case, $\check B_p \< B^*_p$.

It remains to prove that $\check B_p \> B^*_p$. By Lemma 6.6, $a_p+\sqrt p a_2 \< \check B_p (a_p,a_2)$ since $\sqrt c \e ^{p/c} \ge \sqrt p$ if $1\le c\le p$. Thus [in view of (6.22), (6.23) for $q\in I_p$], we only need to show that in the cases (6.17) and (6.18), 
$$ {pa_p \over \ln \( 2+\sqrt p {a_p \over a_2} \) } \< \check B_p (a_p,a_2). \leqno (6.24)$$

If (6.17) is true, then, as it was said above in this proof,  $\check B_p (a_p,a_2) =a_p$; on the other hand, ${a_2 \over a_p} =q <\e^{-p}$, and so, $\ln \( 2+\sqrt p {a_p \over a_2} \) > \ln (2+\sqrt p \e ^p) >p$; now, (6.24) follows.

If (6.18) is true, then $\check B_p (a_p,a_2) \asymp \sqrt p a_2$; on the other hand, ${a_2 \over a_p} =q > \sqrt p / \e$, 
$$ {pa_p \over \ln \( 2+\sqrt p {a_p \over a_2} \) } < {\e  \over \ln2}\, a_2 \,\sqrt p \, ,$$
and, again, (6.24) follows. \sqr

P{\eightrm ROOF OF} T{\eightrm HEOREM}~6.1. Observe that
$$B_p \> B_{p, {ind}, 0} \> \check B_p \> \hat  B_p \> B_p;$$
the first inequality in this chain is trivial; then we successively
apply
Lemma~6.7, Lemma~6.6 and (4.1). Hence,
$$B_p \asymp B_{p,{ind},0} \asymp \check B_p \asymp \hat  B_p.
\leqno
(6.25)$$
Note that for any separable Banach space  $\X $,
$$B_{p, {ind},0} \< B^{\X}_{p, {ind}}\leqno (6.26)$$
since $\Re
%{\Bbb R}
$ may be isometrically embedded into $\X $ and, for
$f\in
{\cal M}  _{ind}(\Re
%{\Bbb R}
)$,
$$\|  \; |f_n+x| - \E|f_n+x| \; \|  _p \to \|  f_n\|  _p \quad\hbox{as}\quad
x\to
\infty.$$
It follows from Theorem 5.1 that
$$B^\X _{p,{ind}} \< \hat  B_p.\leqno (6.27)$$
Now  (6.25)-(6.27) yield
$$B_p \asymp B^\X _{p,{ind}} \asymp B_{p,{ind},0} \asymp \check B_p
\asymp
\hat  B_p.$$
Analogously,
$$B_{p,S}\asymp B^\X _{p, {ind}, S} \asymp B_{p, {ind}, 0, S} \asymp
\check
B_p \asymp \hat  B_p;$$
here, we need to refer to Lemma 6.8 instead of Lemma~6.7 and also to
$\|  d^*\|  _p \le \|S_p\|_p$. It remains to recall Lemma~6.9. \sqr%\square
\medskip

P{\eightrm ROOF OF} T{\eightrm HEOREM}~6.2.  For any given $c\in [1,p]$, choose
$a_p>0$,
$a_2>0$ so that
${a_2/ a_p} = \sqrt c\, \e ^{-p/c}$.
Then $ca_p = \sqrt c\, \e ^{p/c}a_2$. Put $\alpha = (c-1)/p$, so that
$c=p\alpha+1$, $0\le \alpha\le 1$. Then 
$$\eqalign{\check B_p (a_p,a_2) &\ge a_pg_p (\alpha) =
a_pc\e ^{-p\alpha/c} \ge 
\e ^{-1} ca_p\cr
&= (2e)^{-1} (ca_p+ \sqrt c\, \e ^{p/c} a_2),}$$
where $g_p$ is	defined by (6.10). It remains to apply Theorem~6.1
and also
to recall that inequalities (4.2)-(4.4), (5.2)-(5.4), (5.8)-(5.10)
were
obtained by choosing particular values of $c$. \sqr

P{\eightrm ROOF OF} P{\eightrm ROPOSITION}~6.3.  For any given $\alpha\in [0,1]$,
choose
$a_p>0$, $a_2>0$ so that
$${a_2\over a_p} = q_p (\alpha),\leqno (6.28)$$
where $q_p (\alpha)$ is defined by (6.11). Then, putting $c=p\alpha
+1$, we
see that $1\le c\le p+1$ and, in view of (6.10), (6.28),
$$\eqalign{a_pg_p (\alpha) = a_p c\exp\left\{{p\alpha\over 2}
{\alpha-2\over p\alpha+1}\right\}\asymp a_pc\asymp ca_p + \sqrt c\, \e ^{p/c}a_2}$$
since, in consequence of (6.28), (6.11),
$$\sqrt c\, \e ^{p/c} a_2 = ca_p \exp\left\{{p\alpha\over
2(p\alpha+1)}\right\} \asymp ca_p.$$
Thus, if the defined above $c\le p$, we see that
$$\hat  B_p \< a_pg_p (\alpha)\leqno (6.29)$$
If however $c\in (p,p+1]$, then
$$a_pg_p (\alpha)\asymp ca_p+\sqrt c\, \e ^{p/c} a_2 \ge pa_p + \sqrt p \, \e ^{p/(p+1)} a_2 \>
\hat 
B_p,$$
so that (6.29) holds. Now it remains to apply Theorem~6.1. \sqr

\vskip .8cm

\n {\bf 7. Chung type  bounds on moments.}
\ \ Consider
$$\eqalign{B^{{\rm Ch}}_p (a_p,n) := &\sup\{\| \; \|  f_n\| \; \| _p\colon \
(0,f_1,\ldots,
f_n,f_n,f_n,\ldots) \in {\cal M}  (\X ),\cr
&\,\qquad 
n \hbox{ is fixed, }\quad  D\, \|S_p\|_p = a_p,\quad \X \ \hbox{is $(2,D)$-smooth}\}.}$$

\proclaim T HEOREM 7.1. For all $p>2, a_p>0$, $n=1,2,\ldots$, one has 
$$B^{{\rm Ch}}_p (a_p,n)
\asymp \sqrt{p\wedge n}\, n^{p-2\over 2p}a_p.$$

P{\eightrm ROOF.}  Let $(0,f_1,\ldots, f_n,f_n,f_n,\ldots)\in {\cal M}  (\X )$.
Using
 H\"older's  inequality two times, we see that $\|  s_2\|  _p \le
n^{(p-2)/(2p)}\|S_p\|_p$. If $p \le 3n^{(p-2)/p}$, then $p\|  d^*\|  _p \le
p\|S_p\|_p \<
p^{1/2} n^{(p-2)/2p}\|S_p\|_p$, and hence, (4.2) and $p\le 3n$ yield
$$B^{{\rm Ch}}_p (a_p,n) \< \sqrt{p\wedge n}\, n^{p-2\over 2p}a_p.\leqno
(7.1)$$
If now $p>3n^{(p-2)/p}$, then $p>3$, $n<p^{p/(p-2)} \< p$. The inequalities of Minkowsky and H\"older
give
$\| \; \|  f_n\| \; \| _p \le \sum^n_{j=1} \| \; \| d_j\| \; \| _p \le n^{(p-1)/p}\|S_p\|_p$;
this and $p\>n$ show that (7.1) holds. It remains to prove that
$$B^{{\rm Ch}}_p (a_p,n) \> \sqrt{p\wedge n}\; n^{p-2\over 2p}a_p.\leqno
(7.2)$$

Let $d_j$'s be independent, $\P (d_j = \pm u) = 1/2$, $u>0$,
$j=1,\ldots, n$.
Then the multinomial formula yields
$$\|  f_n\|  _{2m} = \left(\sum^{m\wedge n}_{j=1} {n\choose j}
j!\Gamma_{j,m}\right)^{1\over 2m} u \ \asymp \left(\sum^{m\wedge
n}_{j=1}
{n\choose j} j^{2m}\right)^{1\over 2m}u, \leqno (7.3)$$
where $ \Gamma_{j,m}$ is defined by (6.7), $m=1,2,\ldots, $ and (6.19) is
used.
Let $m$ be the integer part of $p/2$. If $m\le n$, then (7.3) shows
that
$$\eqalign{\|  f_n\|  _p \ge \|  f_n\|  _{2m} \> {n\choose m} ^{1\over
2m} mu
\ge \left({n\over m}\right)^{1\over 2} mu= \sqrt m\, n^{(p-2)/2p} \|S_p\|_p,}$$
which implies (7.2) in the case $m\le n$. If, finally, $m>n$, then it
follows from (7.3) that
$$\|  f_n\|  _p \ge \|  f_n\|  _{2m} \> {n\choose
n}^{1\over 2m} nu = n^{(p-1)/p}\|S_p\|_p,$$
so that (7.2) holds. \sqr

Chung (1951), pp.348-349, showed that in the case $f\in
{\cal M}  _{ind}(\Re
%{\Bbb R}
)$, 
the inequality of Marcinkiewicz and Zygmund (1937), p.~87, implies an
estimate like (7.1) but with some $C_0(p)$ depending only on $p$,
instead of $\sqrt{p\wedge n}$. As it was pointed out in
Dharmadhikari,
Fabian and Jogdeo (1968), an analogous result for $f\in {\cal M}  (\Re
%{\Bbb R}
)$
is
implied by the generalization of the Marcienkiewicz-Zygmund
inequality
obtained by Burkholder (1966), Theorem~9 [see also Burkholder (1973),
Theorems~3.2, 15.1, or (9.2) below]; using Proposition~9.1, one can
see
that Theorem~15.1 in Burkholder (1973) in fact gives $C_0(p) = p$; a
constant of the same order $p$ is given in Theorem~3.2 of Burkholder
(1973).

The direct proof due to Dharmadhihari, Fabian and Jogdeo (1968) yields
$C_0(p) = p2^p$. For $f\in {\cal M}  _{ind}(\Re
%{\Bbb R}
)$, Dharmadhikari and
Jogdeo
(1969) obtained $C_0(p) = p/\ln p$.

For $f\in {\cal M}  _{ind}(\Re
%{\Bbb R}
)$, the result of Whittle (1960) implies
the
Marcinkiewicz-Zygmund inequality with the best constant and, along
with
the mentioned above Chung's remark (see also Ros\'en (1970)), leads
to (7.1)
but with $\sqrt p$ instead of $\sqrt{p\wedge n}$; so, for $n>p$, it
gives
the optimum.

What has been said is a reason for the referring to (7.1) as to an
optimum
bound on moments of the Chung type.

\n R{\eightrm EMARK.} Bounds of the Chung type on central moments
of the
norm of the sum of independent random vectors in any separable Banach
space can be easily derived from Theorem~7.1 (cf.
Theorem~5.1).  \sqr
 
\vskip .8cm

\n {\bf 8. One-sided bounds for the distributions of real-valued
(super)martingales.}
\ \ Let ${\cal M}  _-$ stand for the set of all real-valued supermartingales $f\in {\cal S} (\Re)$.  For \break $f\in {\cal M}  _-$, put
$f^*_+ =
\sup\limits_j f_j$, $d^*_+ = \sup\limits_j d_j$.

\proclaim T HEOREM 8.1. If $f\in {\cal M}  _-$ and $\lambda>0$ is such that
$\E\e ^{\lambda d_j} < \infty$ for each $j$, then for all $r\ge 0$,
$$\eqalign{\P (f^*_+\ge  r) &\le \e ^{-\lambda r} \left\| \prod
^\infty_{j=1} (1+e _j)\right\|_\infty\cr
&\le \exp\left\{-\lambda r + \left\| \sum^\infty _{j=1}
e _j\right\|_\infty\right\},}$$
where $e _j := \E _{j-1} (\e ^{\lambda d_j} -1-\lambda d_j)$.

P{\eightrm ROOF}  follows from the trivial remark that the sequence:\ $G_0 := 1$,
\hfil\break $G_j := \e ^{\lambda f_j} \prod^j_{i=1} (1+e _i)^{-1},\ 
j=1,2,\ldots$,  
is a positive supermartingale [cf. the end of the proof of Theorem 3.2]. \sqr 

\proclaim T HEOREM 8.2. Suppose that $f\in {\cal M}  _-$, $\|  d^*_+\|  _\infty
\le a$,
$\|  s_2\|  _\infty \le b$ for some $a>0$, $b>0$. Then for all $r\ge
0$,
$$\P (f^*_+ \ge r) \le \exp \left[{r\over a} - \left({r\over a} +
{b^2\over
a^2}\right) \ln \left(1+{ra\over b^2}\right)\right] \le
\left({\e b^2\over
ra}\right)^{r/a}.$$

P{\eightrm ROOF}  is quite similar to that of Theorem~3.4. \sqr 

\proclaim T HEOREM 8.3. If $f\in {\cal M}  _-$, $p\ge 2$, $1\le c \le p$, then
$$\|  f^*_+\|  _p \< c\|  d^*_+\|  _p + \sqrt c\, \e ^{p/c} \|  s_2\|  _p.$$

%%%P{\eightrm ROOF.}
 We use

\proclaim L EMMA 8.4. (cf. Lemma 4.2.) If $\lambda>0$, $\delta_1>0$,
$\delta_2>0$, $\beta-1-\delta_2>0$, $f\in {\cal M}  _-$, then
$$\P (f^*_+ > \beta\lambda, w^*_+ \le \lambda) \le \varepsilon \P (f^*_+
>\lambda),$$
where
$$\eqalign{w^*_+ &= \left({d^*_+\over \delta_2}\right)\vee
\left({s_2\over
\delta_1}\right),\cr
\vp &= \left({e\over N} {\delta^2_1\over \delta^2_2}\right)^N, \quad N =
{\beta-1-\delta_2\over \delta_2}.}$$

P{\eightrm ROOF.} Put $\bar d_j = d_j I\{d_j \le \delta_2\lambda\}$,
$\bar
f_j = \bar d_0 + \cdots + \bar d_j$, $\bar s_{2,j+1} =
(\sum^{j+1}_{i=1} \E _{i-1} \bar d^2_i)^{1/2}$,
\break $\bar h_j = \bar f_{(j\wedge \tau\wedge\nu)\vee \mu} - \bar
f_\mu,\ 
j=0,1,\ldots\ $, where $\mu = \inf\{j\colon \ \bar f_j>\lambda\},\quad \ $ 
$\nu = \inf\{j\colon \ \bar f_j>\beta\lambda\},\quad \ $ 
$\tau = \inf\{j\colon \ \bar s_{2,j+1} > \delta_1\lambda\}$.
Then $(\bar h_j)$ is a supermartingale conditionally on $F_\mu$, and
$$\eqalign{\P (f^*_+ > \beta\lambda, w^*_+ \le \lambda) &= \P (\bar f^*_+
>
\beta \lambda, w^*_+ \le \lambda)\cr
&\le \P (\bar h^*_+ >(\beta-1-\delta_2)\lambda)\cr
&= \E\,\P (h^*_+ > (\beta -1-\delta_2) \lambda\mid F_\mu)
I\{\mu<\infty\}\cr
&= \vp \P (\bar f^*_+ > \lambda) \le \vp \P (f^*_+>\lambda);}$$
here, we put $\bar f^*_+ = \sup_j \bar f_j$, $\bar h^*_+ = \sup_j \bar h_j$ and took into
account
Theorem~8.2. \sqr 

P{\eightrm ROOF OF} T{\eightrm HEOREM}~8.3 is similar to that
of
Theorem~4.1, but  simpler. Here, we do not need to symmetrize; instead,
we can
apply Lemma~7.1 of Burkholder (1973) directly to $f^*_+$. \sqr

The following is a refinement of Theorem 8.1. 

\proclaim T HEOREM 8.5. Let $f\in {\cal M}  (\Re
%{\Bbb R}
)$, $d_j = u_j -
\E _{j-1}u_j$,
$\lambda>0$, $\E\e ^{\lambda u_j} < \infty$, $j=1,2,\ldots\ $. Then for
all
$r\ge 0$,
$$\P (f^*_+\ge r) \le \exp \left\{ -\lambda r + \left\|
\sum^\infty_{j=1} e _j \right\| _\infty\right\},$$
where $e _j = \E _{j-1}(\e ^{\lambda u_j} - 1-\lambda u_j)$. If, besides,
$\E _{j-1} u_j\ge 0$ for all $j=1,2,\ldots,$ then, moreover,
$$\P (f^*_+\ge r) \le \e ^{-\lambda r} \left\| \prod^\infty_{j=1}
(1+e _j)\right\|_\infty.$$

P{\eightrm ROOF}  is analogous to that in Pinelis and Sakhanenko (1985) but simpler. The elementary inequalities
$$\eqalign{(1+a+b) \e ^{-b} &\le \e ^a \quad (a\in \Re
,\ b\in \Re
),\cr
(1+a+b) \e ^{-b} &\le 1+a \quad (a\ge 0,\ b\ge 0)\cr}$$
imply, resp., that 
$$\E _{j-1}\e ^{\lambda d_j} \le \e ^{e _j},$$
and, if $\E _{j-1}u_j\ge 0$,
$$\E _{j-1} \e ^{\lambda d_j} \le 1+e _j,$$
$j=1,2,\ldots,$ if one chooses $a=e _j$, $b=\lambda \E _{j-1}u_j$.
Hence,
$$\exp\left\{\lambda f_j  - \sum^j_{i=1} e _j\right\},\qquad
j=0,1,2,\ldots,$$
is a supermartingale and, if $\E _{j-1} u_j\ge 0 \quad \forall j$, so is
$$\e ^{\lambda f_j} \prod ^j_{i=1} (1+e _i)^{-1},\qquad j=0,1,2,\ldots \
.$$
It remains to use reasoning like that at the end of the proof of
Theorem~3.2. \sqr

\n R{\eightrm EMARK.} Martingales like those in Theorem~8.5 may arise,
e.g., as
a result of truncating and subsequent centering of the increments of
other
martingales. The aim of Theorem~8.5 is to provide for the best constants
in
 exponential inequalities for martingales, which cannot be
reached
via, e.g., the straightforward estimate $|d_j|\le |u_j| +
|\E _{j-1}u_j|$. \sqr

 As
illustration, let us give the following corollaries to Theorem~8.5,
which
are refinements of Theorems~3.3, 3.4, resp., for the particular case
$\X =\Re
%{\Bbb R}
$.

\proclaim T HEOREM 8.6. Suppose that $f\in {\cal M}  (\Re
%{\Bbb R}
)$, $d_j =
u_j-\E _{j-1}u_j$, $j=1,2,\ldots,$ and
$$\left\|\sum^\infty_{j=1} \E _{j-1} |u_j|^m
\right\|_\infty \le
m! \Gamma^{m-2} B^2/2$$
for some $\Gamma>0$, $B>0$ and $m=2,3,\ldots$\ . Then for all $r\ge
0$, 
$$\P (f^*_+ \ge r) \le \exp\left(- {r^2\over B^2+B\sqrt{B^2+2\Gamma
r}}\right).$$

P{\eightrm ROOF}   is almost literally the same as that of Theorem~3.3, but rests
upon
Theorem~8.5, instead of Theorem 3.1. \sqr

\proclaim T HEOREM 8.7. Suppose that $f\in {\cal M}  (\Re
%{\Bbb R}
)$, $d_j = u_j -
\E _{j-1}u_j$, $j=1,2,\ldots$, and $\|  u^*\|  _\infty \le ~a$,
$\left\|\sum\limits^\infty_{j=1}
\E _{j-1}u^2_j\right\|_\infty
\le b^2$ for some $a>0$, $b>0$. Then for all $r\ge 0$,
$$\eqalign{\P (f^*_+\ge r) &\le \exp \left[{r\over a} - \left({r\over
a} +
{b^2\over a^2}\right) \ln \left( 1+{ra\over b^2}\right)\right]\cr
&\le \left({\e b^2\over ra}\right)^{r/a}.}$$

P{\eightrm ROOF}   differs from that of Theorem~3.4 only in that we use
Theorem~8.5
in place of Theorem~3.1.\sqr

If, instead of Theorem~8.5, we had used the ``naive'' estimate $|d_j|
\le
|u_j| + |\E _{j-1}u_j|$ and, say, Theorem~3.1, we would have hardly been
able
to obtain inequalities better than ones like those in Theorems~8.6,
8.7 but
with $2\Gamma$, $2B$, $2a$, $2b$ in place of $\Gamma, B, a, b$, resp.
The
gain provided, e.g.,  by Theorem 8.6 is quite significant. If, say, 
$$\exp\left(-{r^2\over (2B)^2 + 2B\sqrt{(2B)^2 +2(2\Gamma)r}}\right)
=
10^{-1},$$
then the bound given in Theorem 8.6 varies from $10^{-2\sqrt 2} < (1.5) 10^{-3}$
(when
$\Gamma r$ is much greater than $B^2$) to $10^{-4}$ (when $B^2$ is
much
greater than $\Gamma r$). 

\vskip .8cm

\n {\bf 9. Appendix.}
\ \ Here, we find explicitly the best constants in the mentioned in the 
preceding sections bounds due to Burkholder (1973) 	obtainable via  the methods used therein. Besides, we prove that, at least for the
sums
of independent random vectors in Hilbert spaces, the bound on moments
found
in Pinelis (1980) is equivalent to (4.1).

It follows from Theorems 15.1 and 21.1 of Burkholder (1973) that for
$f\in
{\cal M}  (\Re
%{\Bbb R}
)$, $p\ge 1$,
$$\leqalignno{\|  f^*\|  _p &\< C_1 (\|  d^*\|  _p +
\|  s_2\|  _p),&(9.1)\cr
\|  f^*\|  _p &\< (C_2+p) \left\| S_2 \right\|_p;&(9.2)}$$
analysis of the proof in Burkholder (1973) shows that one can put
$$\eqalign{C_i = \inf \left\{{\beta/\delta\over [1-i\beta^p
\delta^2/(\beta-1-\delta)^2]^{1/p}}\colon \ \right. & \beta-1-\delta>0,\ 
\delta>0,\cr
&\ \left. {  {   i\beta^p\delta^2  }  \over  {   (\beta-1-\delta)^2   }   }
<1\right\}, \quad
i=1,2.}$$

\proclaim P ROPOSITION 9.1. For $i=1,2$,\ \  $p\ge 1$, one has $C_i\asymp
p$.

P{\eightrm ROOF.} Suppose that $\beta-1-\delta>0, \delta>0$,
$i\beta^p\delta^2/(\beta-1-\delta)^2<1$, $i=1$ or $2$. Then
$\beta^{p-2}\delta^2<1$, and we see that $C_1>\beta/\delta
>\beta^{p/2}$.
If $\beta\ge  2$, then $C_1>2^{p/2}\>p$. Assume now that $\beta<2$.
Put
$c=p\delta$, $k=(\beta-1-\delta)/\delta$. Then $\beta=1+(k+1)c/p \ge$
$\e ^{(k+1)c/2p}$ since $1+\alpha >\e ^{\alpha/2}$ for $0 < \alpha <1$.
Hence,
$1 > \beta^p\delta^2/(\beta-1-\delta)^2 =$ $\beta^p/k^2\ge
k^{-2}\e ^{(k+1)c/2}$, $c<4(k+1)^{-1} \ln k <4$, $\delta = c/p <4/p$;
therefore,
$$C_1 >1/\delta > p/4.\leqno (9.3)$$
On the other hand, $\beta\le \e ^{(k+1)c/p}$ since $1+\alpha \le \e ^\alpha$,
and so,
for $k=3$, $c=1/4$, one has $\beta\le e$,
$2\beta^p\delta^2/(\beta-1-\delta)^2 = 2\beta^p/k^2 \le
2\e ^{(k+1)c}/k^2<2/3$. Hence,
$$C_2 \le 3\beta/\delta < 3\e/\delta = 3\e p/c =12\e p.\leqno (9.4)$$
Note that $C_1\le C_2$. Now (9.3), (9.4) imply the proposition. \sqr

Finally, we give an interpretation of results of Pinelis (1980). Let
${\cal
F}(p)$, $p>0$ stand for the set of all the pairs $(c_1,c_2)$ that may
be
represented in the form
$$\leqalignno{c_1 &= \sum^{m-1}_{i=0} { c_i(p) b^{-2i}_i \over i!},&(9.5)\cr
c_2 &= {\tilde c_m(p)\over m!} + \sum^{m-1}_{i=0} { c_i(p)
b^{p-2i-2}_i \over i! } ,&(9.6)}$$
where $m$ is the integer part of $p/2$, $b_i>0$,
$$\leqalignno{c_i(p) &= v(p-2i)2^{-i}(2i)! {p\choose 2i}
\prod^{i-1}_{k=0} u(p-2k),&(9.7)\cr
\tilde c_m(p) &= 2^{-m}(2m)! {p\choose 2m} \prod^{m-1}_{i=0}
u(p-2i), %\quad 
&(9.8) \cr
{p\choose k} &= {p (p-1)\ldots (p-(k-1))\over k!},%&(9.8)
}$$
and $u(\cdot)$, $v(\cdot)$ are any functions such that
$$\leqalignno{&u(2) + v(2) \ge 1,\quad  u(2) > 0, \quad v(2) >0,&(9.9)\cr
&u(s)\ge 1, \quad v(s)\ge 1\quad \hbox{for}\quad s>2,&(9.10)\cr
&u(s)^{1\over 3-s} + v(s)^{1\over 3-s} \le 1\quad \hbox{for}\quad
s>3.&(9.11)}$$
 For $a_p>0, a_2>0, p>0$, put
$$B_{p,1980}(a_p,a_2) = \inf\{(c_1a^p_p +c_2a^p_2)^{1/p}\colon \
(c_1,c_2)
\in {\cal F}(p)\}.$$

In Pinelis (1980), bounds on moments of martingales in Hilbert spaces
with
bounded conditional second moments were found. If applied, in
particular,
to any $f\in {\cal M}  _{ind}(\X )$, where $\X $ is any separable Hilbert space,
the
Corollary in Pinelis (1980) states that
$$\|  \; \| f_n\| \; \|  _p \le B_{p,1980} (\|S_p\|_p,\|S_2\|_2),\quad p>0.$$
In turn, this implies (cf. Theorems~5.1, 6.1)
$$B_{p, {ind}, S}(a_p,a_2) \le 2\, B_{p,1980}(a_p,a_2) \leqno (9.12)$$
for any $a_p>0, a_2>0, p>0$ and any separable Banach space $\X $.

\proclaim P ROPOSITION 9.2. For $a_p>0, a_2>0,p>2$,
$$B_{p,1980}(a_p,a_2) \asymp B_p (a_p,a_2).$$
In view of Theorem 6.1, any of the following functions:\ $B_p^*$, $B^\X _{p,
{ind}},
B_{p, {ind}, 0}$, $B_{p,S}$, $B^\X _{p, {ind}, S}$, $B_{p, {ind}, 0,
S}$,
$\check B_p$, $\hat B_p$ may be used here in place of $B_p$.

P{\eightrm ROOF.} For $s\ge 0$, put
$$y_s = g_s(\alpha_s) = \max\{g_s(\alpha)\colon \ 0 \le \alpha\le
1\},\qquad 0 \le \alpha_s \le 1,\leqno (9.13)$$
and, for $s\ge 2$,
$$u(s) = {  { \e ^3y^s_s  } \over  { sq^2 y^{s-2}_{s-2}  }  },\quad v(s) = 10^s y^{s-2}_s,
\leqno
(9.14)$$
where (cf. (6.10))
$$g_s(\alpha) = (s\alpha +1)^{1-\alpha/2}q^\alpha,\quad  q=a_2/a_p.\leqno
(9.15)$$
Note that for $s\ge 0$,
$$y_s\ge g_s(0) \vee g_s(1) = 1\vee (s+1)^{1/2}q.\leqno (9.16)$$
Besides, $y_s$ is increasing in $s$ because so is $g_s(\alpha)$.
Hence, for
$s\ge 2$, one has $u(s)\ge \e ^3y^2_s/(sq^2)\ge \e ^3>1$, $v(s) \le
2\cdot 10^s
y^s_s/(1+q^2)$, $v(s) \ge (10^{s-3}y^{s-3}_s)\vee 1$, and so, for
$s>3$,
$$u(s)^{1\over 3-s} + v(s)^{1\over 3-s} \le \left[{  {sq^2} \over {\e ^3y^2_s  }  } \right]^{1\over
s} +
{1\over 10y_s}\; .\leqno (9.17)$$

Consider first the case $q\in I_s$, where $I_s$ is defined by (6.13).
Then,
$y_s = g_s(\alpha_s) > \e ^{-1}(s\alpha_s +1)$
[cf. (6.15)], and, in view of (6.11),
$$\eqalign{q=q_s(\alpha_s) < \sqrt{s\alpha_s+1}\, \exp\left\{{1\over
2} -
{s\over s\alpha_s+1}\right\}< \sqrt{y_s} \e ^{1-s/(y_s\e)};}$$
hence,
$$\left({sq^2\over \e ^3y^2_s}\right)^{1\over s} < \left({s\over
\e y_s}\right)^{1\over s} \e ^{-1/(y_s\e)} \le \exp\left\{ {1\over \e ^2y_s}
-
{1\over \e y_s}\right\} < 1- {1\over 10y_s}$$
since $(s/a)^{1/s} \le \e ^{1/(a\e)}$ for $a>0$ and $\e ^{-\alpha} < 1 -
\alpha/2$ for $0 < \alpha < 1$. Thus, (9.17) yields (9.11) if
$q\in
I_s$.

If $q< \min  I_s = \e ^{-s}$, then $y_s=g_s(0) = 1$, and
$$\left({sq^2\over \e ^3y^2_s}\right)^{1\over s} < \left({s\over
\e ^3}\right)^{1\over s} \e ^{-2} \le \e ^{\e ^{-4}-2} < \e ^{-1},$$
so that (9.17) again yields (9.11).

Now suppose that $q>\max I_s(> [(s+1)/\e]^{1/2})$. Then $y_s = g_s(1)
=
(s+1)^{1/2} q> s^{1/2}q$, and for $s> 3$,\ \  $[(sq^2)/
(\e ^3y^2_s)]^{1\over s} < \e ^{-3/s} < 1  - 1/ s$,
$${1\over 10y_s} = {1\over 10\sqrt{s+1}\, q} < {\sqrt e\over 10(s+1)}
<
{1\over s},$$
so that  (9.17) once more yields (9.11).

Thus, $u(s)$ and $v(s)$ defined by (9.14) satisfy (9.9)-(9.11).

Let us
write ${\cal E}_1 \ll {\cal E}_2$ if ${\cal E}^{1/p}_1 \< {\cal
E}^{1/p}_2$,
$p>2$.

Using (9.7), (9.14), (9.16) and the Stirling's formula, we see
that for
$p>2$,\break $i=1,\ldots, m-1$,
$$\leqalignno{{c_i(p)\over i!} &\ll y^p_p y^{-2}_{p-2i} q^{-2i} \ll
y^p_p (1+q)^{-2} q^{-2i}&(9.18)\cr
&\le y^p_p q^{-2pi/(p-2)};}$$
similarly, (9.8), (9.14), (9.16) imply
$${\tilde c_m(p)\over m!} \ll y^p_p q^{-2m} y^{-(p-2m)}_{p-2m} \le
y^p_p
q^{-p}.\leqno (9.19)$$
Choose now $b_i \equiv q^{-p/(p-2)}$ so that $b^{-2i}_i a^p_p
=b^{p-2i-2}_ia^p_2$,
$i=1,\ldots, m$, in accordance with (9.15); then (9.5), (9.6),
(9.18),
(9.19), (9.13), (6.12) yield
$(c_1a^p_p + c_2a^p_2)^{1/p} \< y_pa_p = \check B_p$.
Thus, $B_{p,1980}(a_p,a_2) \< \check B_p (a_p,a_2)$. It remains to
recall
(9.12) and Theorem~6.1. \sqr

\vskip .8cm

\vfil\eject

\def\item#1{\hangindent=2cm \hangafter=1}

\centerline{REFERENCES}

\item{1.} Bennett, G. (1962). Probability inequalities for the sum of
independent random variables. {\it J. Amer. Statist. Assoc. \/}{\bf 57}
33-45.

\item{2.} Berger, E. (1991). Marjorization, exponential inequalities
and
almost sure behavior of vector-valued random variables. {\it Ann. Probab.\/}
{\bf
19} 1206-1226.

%%%%\item{3.} Bernstein, S. (1924). {\it Sur une modification de l'
%%%in\'equalit\'e
%%%%de Tchebichef.\/} Annals Science Institute Sov. Ukraine, Sect. Math. I.

\item{4.} Billingsley, P. (1968). Convergence of probability
measures. John
Wiley, New York.

\item{5.} Burkholder, D. (1966). Martingale transforms. {\it  Ann. Math.
Statist.\/} {\bf 37} 1494-1504.

\item{6.} Burkholder, D. (1973). Distribution function inequalities
for
martingales. {\it Ann. Probab.\/} {\bf 1} 19-42.

\item{6a.} Burkholder, D. (1991). Explorations in martingale theory and its applications. Saint-Flour Lectures (1989). {\it Lecture Notes in Mathematics\/} {\bf 1464}\ \ 1-66.

\item{7.} Chung, K.L. (1951). The strong law of large numbers. {\it Proc.
Second
Berkeley Symp. Math. Statist. Prob.\/} 341-352. Univ. of California
Press.

\item{8.} Davis, B. (1970). On the integrability of the martingale
square
function. {\it Israel J. Math.\/} {\bf 8} 187-190.

\item{9.} De Acosta, A. (1981). Inequalities for $B$-valued random
vectors
with applications to the strong law of large numbers. {\it Ann. Probab.\/}
{\bf 9},
157-161.

\item{10.} Dharmadhikari, S.W., Fabian, V. and Jogdeo, K. (1968).
Bounds on
the moments of martingales. {\it Ann. Math. Statist.\/} {\bf 39} 1719-1723.

\item{11.} Dharmadhikari, S.W. and Jogdeo, K. (1969). Bounds on
moments of
certain random variables. {\it Ann. Math. Statist.\/} {\bf 40} 1506-1508.

\item {\ } Haagerup, U.  (1982). The best constants in the Khintchine inequality. {\it Studia Math.\/} {\bf 70}\ 232-283.

\item{12.} Hall, P. and Marron, J.S. (1987). Choice of kernel order
in
density estimation. {\it Ann. Statist.\/} {\bf 16} 161-173.

\item{13.} Hitczenko, P. (1990). Best constants in martingale version
of
Rosenthal's inequality. {\it Ann. Probab.\/} {\bf 18} 1656-1668.

\item{13.} Hitczenko, P. (1990a). 
Upper bounds for the $L_p$-norms of
martingales. {\it Probab.  Th.  Rel.  Fields\/}   {\bf 86}  225-237.

\item{14.} Hitczenko, P. (1991). Private communication.

\item{15.}  Hoeffding, W. (1963). Probability inequalities for sums
of
bounded random variables. {\it J. Amer. Statist. Assoc.\/} {\bf 58} 13-30.

\item{16.} Hoffmann-J\o rgensen, J. (1974). Sums of independent
Banach
space valued random variables. {\it Studia Math.\/} {\bf 52} 159-186.

\item{17.} Johnson, W.B., Schechtman, G. and Zinn, J. (1985). Best
constants in moment inequalities for linear combination of
independent and
exchangeable random variables. {\it Ann. Probab.\/} {\bf 13} 234-253.

%%\item{17a.} Kahane, J.-P. (1968). Some random series of functions.
%%Heath,
%%Lexington, Mass.

\item{17b.} Kallenberg, O. and Sztencel, R. (1991). Some
dimension-free
features of vector-valued martingales. {\it Probab. Th. Rel. Fields\/} {\bf
88} 215-247.

\item{18.} Kwapie\'n, S. and Szulga, J. (1991). Hypercontraction
methods in
moment inequalities for series of independent random variables in
normed
spaces. {\it Ann. Probab.\/} {\bf 19} 369-379.

\item{19.} Kwapie\'n, S. and Woyczy\'nski, W.A. (1989). Tangent
sequences
of random variables:\ Basic inequalities and their applications.
{\it Proc.
Conf. on Almost Everywhere Convergence in Probab. and Ergodic Theory,\/}
237-265. Academic, New York.

\item{20.} Marcinkiewicz, J. and Zygmund, A. (1937), Sur les
fonctions
independantes. {\it Fund. Math.\/} {\bf 89} 60-90.

\item{21.} Nagaev, S.V. and Pinelis, I.F. (1977). Some inequalities
for the
distributions of sums of independent random variables. {\it Theory Probab.
Appl.\/}
{\bf 22} 248-256.

\item{22.} Pinelis, I.F. (1974). On some inequalities for the
distributions
of sums of independent random variables. {\it Diploma work.\/} Novosibirsk
State
University. (Russian).

%%\vfill
%%\eject

\item{23.} Pinelis, I.F. (1978). Distributions of sums of independent
random variables with values in a Banach space. {\it Theory Probab. Appl.\/}
{\bf
23} 608-615.

\item{24.} Pinelis, I.F. (1980) Estimates of moments of
infinite-dimensional martingales. {\it Math. Notes\/} {\bf 27} 459-462.

\item{25.} Pinelis, I.F. (1981). Limit theorems on large deviations
for
sums of infinite-dimensional random variables when the Cram\'er's
condition
is violated. {\it Manuscript deposited at All-Union Institute for
Scientific and
Technical Information \break (VINITI),\/} no. 1674-81. (Russian).

\item{26.} Pinelis, I.F. (1986). Probability inequalities for sums of
independent random variables with values in a Banach space. {\it Math.
Notes\/} 
{\bf 39} 241-244.

\item {}  Pinelis, I.F. (1990). Inequalities for the distributions of sums of independent random vectors and their applications to
the density estimation. {\it Theory Probab. Appl.\/} {\bf 35}  500-512.

\item{27.}  Pinelis, I. (1992). An approach to inequalities for the distributions of \hfil\break infinite-dimensional martingales. {\it Probability in Banach Spaces, {\bf 8}: \ Proceedings of the  Eighth   International  Conference.}\ 128-134. Birkh\"auser, Boston.

\item{28.} Pinelis, I. (1992a). On a distribution inequality for sums of independent  random vectors. {\it Submitted.\/}

\item{29.} Pinelis, I.F. and Sakhanenko, A.I. (1985). Remarks on
inequalities for large deviation probabilities. {\it Theory Probab. Appl.\/}
{\bf
30} 143-148.

\item{30.} Pinelis, I.F. and Utev, S.A. (1984). Estimates of the
moments of
sums of independent random variables. {\it Theory Probab. Appl.\/} {\bf 29}
574-577.

\item{31.} Pinelis, I.F. and Utev, S.A. (1989). Exact exponential
bounds
for sums of independent random variables. {\it Theory Probab. Appl.\/} {\bf
34}
340-346.

\item {\ } Pisier, G. (1975). Martingales with values in uniformly convex spaces. {\it  Israel J. Math.\/} {\bf 20}\ 326-350.

\item{32.} Ros\'en, B. (1970). On bounds on the central moments of
even
order of a sum of independent random variables. {\it Ann. Math. Statist.\/}
{\bf
41} 1074-1077.

\item{33.} Rosenthal, H.P. (1970). On the subspaces of $L^p (p>2)$
spanned
by sequences of independent random variables. {\it Israel J. Math.\/} {\bf 8}
273-303.

\item{34.} Sazonov, V.V. (1974). On the estimation of moments of sums
of
independent random variables. {\it Theory Probab. Appl.\/} {\bf 19} 371-374.

\item{35.} Talagrand, M. (1989). Isoperimetric and integrability of
the sum
of independent Banach-space valued random variables. {\it Ann. Probab.\/}
{\bf 17}
1546-1570.

\item{36.} Wang, M. (1991). Sharp square function inequlities for conditionally symmetric martngales. {\it Trans. AMS\/}
{\bf 328}
393-419.

\item{37.} Whittle, P. (1960). Bounds for the moments of linear and
quadratic forms in independent random variables. {\it Theory Probab. Appl.\/}
{\bf
5} 302-305.

\item{38.} Yurinskii, V.V. (1974).  Exponential bounds for large
deviations. {\it Theory Probab. Appl.\/} {\bf 19} 154-155.

\vskip 1cm

\parskip=0pt
\parindent=7cm \def\\{\eightrm}
D{\\EPARTMENT} {\\OF} M{\\ATHEMATICAL} S{\\CIENCES}

M{\\ICHIGAN} T{\\ECHNOLOGICAL} U{\\NIVERSITY}

H{\\OUGHTON,} M{\\ICHIGAN} 49931

\bye
\end{plain}
\end{document}